\documentclass[journal]{new-aiaa}

\usepackage{xcolor}
\usepackage{hyperref}
\usepackage{tabulary}
\usepackage{graphicx}
\usepackage{mathtools}
\usepackage{cleveref}
\usepackage{amsmath}
\usepackage{amsfonts}
\usepackage{xcolor}
\usepackage{algorithm}
\usepackage{algorithmic}
\usepackage{mathtools} 
\usepackage{tikz}
\usetikzlibrary{positioning, arrows.meta, shapes.geometric, fit}

\algsetup{indent=.5em}

\newcommand{\vect}[1]{\boldsymbol{#1}}
\usepackage{xparse}
\usepackage{etoolbox}

\NewDocumentCommand{\mylist}{m}{%
  \ifstrequal{#1}{\Theta}
    {\vect{\Theta}}
    {\mathcal{\MakeUppercase{#1}}}
}

\newcommand\scalemath[2]{\scalebox{#1}{\mbox{\ensuremath{\displaystyle #2}}}}
\newcommand{\algfunc}[1]{\scalemath{.9}{\mathtt{#1}}}

\newcommand{\red}[1]{\textcolor{black}{#1}}

\title{Cooperative Multi-Agent Path Planning for Heterogeneous UAVs in Contested Environments}
\author{Grant Stagg \footnote{Graduate Student, Department of Electrical and Computer Engineering, {\tt\small ggs24@byu.edu}} and Cameron K. Peterson\footnote{Associate Professor, Department of  Electrical and Computer Engineering,
{\tt\small cammy.peterson@byu.edu}}}
\affil{Brigham Young University, Provo, Utah, 84602}

\begin{document}

\maketitle

\begin{abstract}
This paper addresses the challenge of navigating unmanned aerial vehicles in contested environments by introducing a cooperative multi-agent framework that increases the likelihood of safe UAV traversal. The approach involves two types of UAVs: low-priority agents that explore and localize threats, and a high-priority agent that navigates safely to its target destination while minimizing the risk of detection by enemy radar systems. The low-priority agents employ a decentralized optimization algorithm to balance exploration, radar localization, and safe path identification for the high-priority agent. For the high-priority agent, two path-planning methods are proposed: one for deterministic scenarios using weighted Voronoi diagrams, and another for uncertain scenarios that leverages generalized Voronoi diagrams (incorporating a non-Euclidean criterion derived from uncertainty in the radar’s probability of detection) alongside probabilistic constraints. Both methods employ optimization techniques to refine the trajectories while accounting for kinematic constraints and radar detection probabilities. Numerical simulations demonstrate the effectiveness of our framework. This research advances UAV path planning methodologies by combining heterogeneous multi-agent cooperation, probabilistic modeling, and optimization to enhance mission success in adversarial environments.
\end{abstract}

\section{Introduction}
Unmanned aerial vehicles (UAVs) are being increasingly deployed in complex and contested environments for missions that include reconnaissance, surveillance, and combat operations. A key challenge for the UAVs in these scenarios is ensuring their safe navigation while avoiding detection by enemy radar systems.
Radar detection is inherently probabilistic, influenced by factors such as radar power, environmental conditions, and UAV positioning. This uncertainty poses a significant challenge for mission-critical UAVs that must traverse hostile regions while minimizing detection risks. 

\begin{figure*}
\includegraphics[width=0.95\linewidth,trim={0.0cm .3cm 0.5cm .1cm},clip]{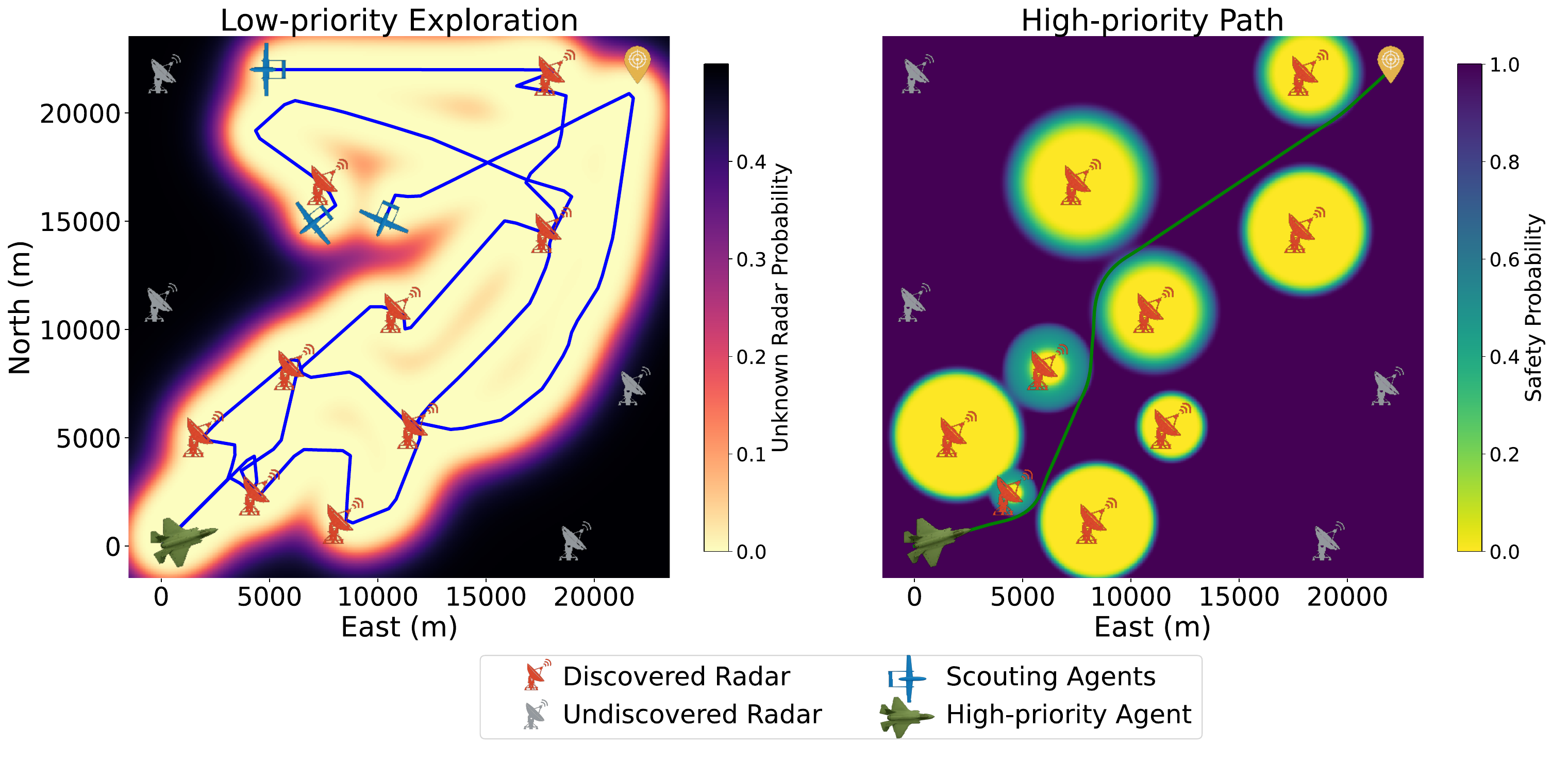}
    \caption{This figure illustrates our low- and high-priority agent framework. Low-priority agents (green) explore the environment, identify radar stations (red), and map out safe and unsafe regions. The high-priority agent (blue) leverages this information to plan a safe path to its goal location (green home icon).}\label{fig:overview}
\end{figure*}

In this paper, we address this challenge by introducing a cooperative framework involving two classes of UAVs: a high-priority agent with a critical mission objective and multiple low-priority agents tasked with scouting the area  (see Figure~\ref{fig:overview}). The scout agents explore the environment, intercept radar emissions, and estimate radar locations and capabilities, providing critical information for the high-priority UAV. The high-priority UAV then leverages this information to traverse a safe path from its initial location to a target destination. 

The use of heterogeneous vehicles increases the likelihood of mission success by designating low-priority vehicles as expendable agents. To maximize their value, we optimize their paths using a three-part objective function that prioritizes uncertainty reduction and exploration while ensuring efficient information gathering along the shortest route. These vehicles integrate probabilistic radar modeling with cooperative multi-agent path planning strategies to gather information that increases the survivability and operational effectiveness of high-priority UAVs operating in adversarial environments.

In addition to our low-priority path planning algorithm, we introduce two high-priority path planning algorithms tailored to different levels of radar parameter certainty. The first algorithm addresses the deterministic case, where all radar parameters are known. First, we use a multiplicatively weighted Voronoi diagram together with an A* search to produce an initial feasible path. That path is then refined via an interior-point optimizer to yield a minimum-time trajectory that meets both safety and kinematic constraints. The second algorithm is designed for scenarios with uncertain radar parameters. In this case, we use a generalized Voronoi diagram and the A* algorithm to find an initial feasible trajectory, followed by an optimization process to refine the path. 
In the generalized Voronoi diagram, each point in the operational domain is associated with the radar system that maximizes the likelihood that the probability of detection exceeds the prescribed threshold.

The key contributions of our work are summarized as follows:
\begin{itemize}
    \item \textbf{Cooperative Path Planning for Low-Priority Agents: } We develop an algorithm that optimally balances regional exploration and the reduction of uncertainty in detected radar locations.
    \item \textbf{Path Planning with Known Radar Parameters: } Utilizing weighted Voronoi diagrams, we design a path planning algorithm that minimizes radar detection risk when radar parameters are fully known.
    \item \textbf{Path Planning with Uncertain Radar Parameters:} For scenarios with uncertain radar parameters, we introduce a path planning algorithm based on generalized Voronoi diagrams to enhance the stealth and survivability of a high-value UAV.
\end{itemize}

This paper is organized as follows. Section~\ref{sec:related_works} reviews relevant past research. In Section~\ref{sec:Problem_statement}, we define the problem and its path-planning applications. Section~\ref{sec:background} provides essential background information on radar localization, detection uncertainty quantification, Voronoi diagrams, and B-splines. Section~\ref{sec:lp-path-planning} introduces our multi-agent path planning algorithm for low-priority vehicles, focusing on exploration and uncertainty reduction. In Section~\ref{sec:HP_path_planning}, we present our high-priority path-planning algorithm for radar detection avoidance in both deterministic and uncertain scenarios. Simulation results are presented in Section~\ref{sec:results}, and finally, Section~\ref{sec:conclusion} concludes the paper.

\section{Related Works}\label{sec:related_works}
Our approach involves two different forms of path planning: one tailored for low-priority attritable agents and the other for a high-priority agent. The design of our low-priority path-planning algorithm is informed by previous research in emitter localization~\cite{ponda2009FIMObjectiveFunctions,xu2024systematical,xu2020rssAndAOA,sahu2022optimal,shahidian2017single,shahidian2015autonomous,dogancay2012uav}. Our high-priority path-planning algorithm draws on previous studies in radar-avoidance path planning~\cite{hameed2022reinforcement,li2022genetic,zhang2020improcedAstar,basbous20182dSimulatedAnnealing,zhao2016fastAstar,gao2014AntColonyOpt,kabamba2006optimal,pelosi2012range,zabarankin2006calculusOfVariation,costley2022path,inanc2008framework,sun2017two} and obstacle-avoidance techniques using Voronoi diagrams~\cite{al2020voronoi,judd2001spline,liu2022fusion,choi2010real,zhang2018quantitative,chen2016path,chen2013research}.

Our approach involves two different forms of path planning: one tailored for low-priority attritible agents, and the other for a high-priority agent. The design of our low-priority path planning algorithm is informed by previous research in emitter localization~\cite{ponda2009FIMObjectiveFunctions,xu2024systematical,xu2020rssAndAOA,sahu2022optimal,shahidian2017single,shahidian2015autonomous,dogancay2012uav}. Our high-priority path-planning algorithm draws on previous studies in radar avoidance path planning~\cite{hameed2022reinforcement,li2022genetic,zhang2020improcedAstar,basbous20182dSimulatedAnnealing,zhao2016fastAstar,gao2014AntColonyOpt,kabamba2006optimal,pelosi2012range,zabarankin2006calculusOfVariation,costley2022path,sun2017two} and obstacle avoidance techniques using Voronoi diagrams~\cite{al2020voronoi,judd2001spline,liu2022fusion,choi2010real,zhang2018quantitative,chen2016path,chen2013research}.

\red{The objectives of the low-priority agents are to explore the environment, discover new radar stations, and reduce uncertainty in the estimated locations of the identified radars.} Low-priority UAVs discover and localize enemy radar stations. 
Prior research in this area employed the extended Kalman filter (EKF) to track targets and radio frequency sources based on the angle of arrival and signal strength  measurements~\cite{xu2024systematical,shahidian2017single,shahidian2015autonomous,ponda2009FIMObjectiveFunctions}. Similarly, we utilize the EKF to estimate the locations of enemy radar stations and their effective radiated power. To aid in localization, low-priority agents travel to measurement locations that give the best information about the environment.

\red{Optimization techniques for measurement placement and trajectory planning are integral to these approaches, often leveraging either the Fisher information matrix (FIM)~\cite{shahidian2015autonomous, xu2024systematical,ponda2009FIMObjectiveFunctions} or the EKF covariance matrix~\cite{sahu2022optimal,tzoreff2017path,ponda2009FIMObjectiveFunctions,xu2020rssAndAOA,shahidian2017single} in their objective functions.} Since scalar representations are required for optimization, common FIM-based scalarizations include the determinant of FIM (D-optimality), the trace of FIM (A-optimality), and the largest eigenvalue of FIM (E-optimality)~\cite{ponda2009FIMObjectiveFunctions}. For EKF-based methods, scalarizations like the determinant~\cite{xu2024systematical} or trace~\cite{shahidian2015autonomous} of the covariance matrix are frequently used in the objective functions. These works aim to optimize path planning by minimizing the scalarization of the EKF covariance matrix or maximizing the scalarization of the FIM. The primary goal is to reduce uncertainty in the estimates as efficiently as possible. However, these works assume there is already an estimate of where the emitter is located, which limits their applicability in scenarios where emitters are initially unknown. Our approach addresses this limitation by incorporating exploration strategies that allow agents to discover and localize previously unidentified radar sources, thus broadening the scope of potential applications.  

Building upon these works, we utilize the reduction of the determinant of the EKF covariance in our objective function. We also include extra terms in the objective that incentivize exploration, creating a trade-off between reducing the uncertainty in the estimate of the already discovered radar and exploring to discover undiscovered radar. We do this through a Bayesian method outlined in Section~\ref{sec:LPP_objective}. We also add a third term in the objective function that promotes exploration towards the goal location of the high-priority agent, facilitating the identification of a safe trajectory for the high-priority agent. Additionally, this term helps prioritize efficient navigation by balancing exploration and safety requirements.


For high-priority agents, our research builds on prior work in path planning using Voronoi diagrams and radar detection avoidance strategies.  The authors in~\cite{al2020voronoi,judd2001spline,liu2022fusion} use Voronoi diagrams to generate initial trajectories, then use smoothing or optimization algorithms to find feasible minimum-time trajectories. \red{Similarly, we employ Voronoi diagrams to generate an initial feasible trajectory and then refine the path using a smoothing and optimization algorithm.}

Past research~\cite{choi2010real,zhang2018quantitative,chen2013research,chen2016path,dai2010path} has used standard Voronoi diagrams to plan the path around the threat points, such as radar, missile, or terrain.  In these studies, threat levels were modeled by assigning weights to edges based on the type and parameters of each threat. The minimum threat path was then computed through the Voronoi diagram. However, these approaches do not account for varying threat levels, e.g., the ideal path should pass closer to a lesser threat. \red{We address this limitation by employing weighted Voronoi diagrams, where the weighted Voronoi ridge naturally shifts closer to the radar with a smaller range, reflecting their reduced threat levels.} Weighted Voronoi diagrams are defined using a multiplicatively weighted distance metric instead of the standard Euclidean distance.

The approach in~\cite{liu2022fusion} utilizes weighted Voronoi diagrams to generate trajectories around obstacles. \red{Extending this, we apply both weighted and generalized (non-Euclidean criterion) Voronoi diagrams to design initial trajectories that avoid radar.}  Instead of navigating around physical obstacles, we leverage weighted Voronoi diagrams to identify the ridge of minimum probability of detection (PD) between two radar stations with different capabilities. To our knowledge, this represents the first application of weighted Voronoi diagrams in this context. 

Prior research in radar detection avoidance path planning has predominantly focused on planning routes through enemy radar under the assumption that the radar parameters, such as location and power, are fully known~\cite{hameed2022reinforcement,li2022genetic,zhang2020improcedAstar,basbous20182dSimulatedAnnealing,zhao2016fastAstar,gao2014AntColonyOpt,kabamba2006optimal,pelosi2012range,zabarankin2006calculusOfVariation}. The authors in these works use similar formulations and minimize radar detection or detection probability using the known information. The methods vary significantly in their approaches. For instance, the authors in~\cite{hameed2022reinforcement} train a deep reinforcement learning model to generate paths to minimize or avoid radar tracking. In~\cite{li2022genetic}, a genetic optimization algorithm is employed to plan safe paths, while leveraging terrain masking to reduce detection risks. The authors in~\cite{zhang2020improcedAstar} and~\cite{zhao2016fastAstar} developed improved A* algorithms to find paths through enemy radar, while simulated annealing optimization is applied in~\cite{basbous20182dSimulatedAnnealing}. Other techniques include ant colony optimization~\cite{gao2014AntColonyOpt} and calculus of variations for optimal path planning~\cite{zabarankin2006calculusOfVariation,kabamba2006optimal}. \red{Despite their methodological diversity, all these approaches share a critical limitation: they assume perfect knowledge of radar parameters.} This assumption neglects the inherent uncertainty in contested environments, where radar locations and capabilities are often unknown or probabilistic.

Recently, the authors in~\cite{costley2022path,costley2022sensitivity,costley2023sensitivity} proposed linearization-based methods to approximate the impact of uncertainty in radar detection models, which share some similarities with the approach presented here. In~\cite{costley2022sensitivity}, they analyze the sensitivity of the radar PD equations with respect to the agent's state, enabling approximation of detection uncertainty arising from variability in the agent trajectory. This analysis is extended in~\cite{costley2023sensitivity} to consider sensitivities with respect to radar parameters, thereby capturing the effect of parameter uncertainty on PD.
Our method builds on these ideas by simultaneously accounting for three classes of variables: (1) known agent states with associated uncertainty (e.g., from navigation or sensing error), (2) estimated radar parameters with uncertainty derived from an estimator, and (3) fully unknown (unobservable) radar states modeled via prior belief distributions.

The path planning approach in~\cite{costley2022path} employs a visibility graph around radar-based avoidance polygons, with smoothing applied to generate a flyable trajectory. Because the visibility graph may place the initial path near safety boundaries, smoothing can occasionally introduce constraint violations, which motivated their iterative approach of expanding polygons and replanning.
In contrast, we utilize Voronoi diagrams to generate an initial feasible trajectory (satisfying kinematic feasibility and maximum PD constraints), followed by refinement using an interior point optimization algorithm. Both visibility-graph and Voronoi-based paths require smoothing to ensure kinematic feasibility. However, the Voronoi-based initialization is naturally biased away from high-risk areas, which reduces the likelihood that smoothing will introduce violations. This, in turn, allows the trajectory to be directly refined with a continuous optimization method, without requiring iterative replanning.

\section{Problem Statement}\label{sec:Problem_statement}
We consider a scenario where all the agents (UAVs) operate in a 2D region $D\subset \mathbb{R}^2$. The low-priority agents are located at $\vect{x}_{l,i} = [x_{l,i},\;y_{l,i}]^\top \forall i \in \{1,\hdots, N_l\}$, where $N_l$ is the number of low-priority agents, $x_{l,i}$ is the distance East of the origin,  and $y_{l,i}$ is the distance North of the origin. Additionally, we assume a single high-priority agent located at $\vect{x}_{h} = [x_{h},\;y_{h}]^\top$. The high-priority agent's objective is to plan a path starting at its initial location $\vect{x}_{h_0}$ and ending at its goal location $\vect{x}_{h_f}$, while avoiding detection by enemy radar. Meanwhile, the low-priority agent's aim is to explore the region, detect and locate enemy radar, and identify a safe path for the high-priority agent.

There are $N_r$ enemy radar stations in $D$ at locations, $\vect{x}_{r,j} = [x_{r,j},y_{r,j}]^\top\;\forall j \in \{1, \hdots, N_r\}$. These locations are unknown to the agents.  The low-priority agents are equipped with sensors capable of intercepting enemy radar emissions. Each agent is assumed to measure the angle of arrival $\phi_{i,j}$ of the enemy radar signals, along with the received power. The received power at the $i^\text{th}$ agent from the $j^\text{th}$ radar is given by
\begin{equation}
    S_{E,i,j} = \frac{P_{T,j}G_{T,j}G_{I ,i}\lambda^2}{(4\pi)^2R_{i,j}^2L_j},
\end{equation}
where $P_{T,j}$ is the power transmitted from the $j^\text{th}$ enemy radar, $G_{T,j}$ is the transmit gain of the $j^\text{th}$ enemy radar, $G_{I,i}$ is the gain of the $i^\text{th}$ agent's intercept antenna, $\lambda$ is the radar wavelength (we assume the wavelength for all radar is the same), $R_{i,j}$ is the distance between the $i^\text{th}$ agent and the $j^\text{th}$ radar, and $L_{j}$ is the path loss. 

We wish to estimate the location of the radar  $\vect{x}_{r,j}$ and the effective radiated power (ERP) $P_{E,j}$ for the $j^{\text{th}}$ radar, where the effective radiated power is
\begin{equation}\label{eq:effective_radiated_power}
    P_{E,j}=\frac{P_{T,j}G_{T,j}}{L_j}.
\end{equation}
We estimate the ERP directly because the transmitted power and antenna gain are not separately identifiable from received measurements.
We combine the location and ERP for each radar into an augmented state that will be estimated. The augmented state for the $j^{\text{th}}$ radar station is denoted by 
$\vect{\theta}_{e,j} = 
[{x}_{r,j}, \;
{y}_{r,j}, \;
{P}_{E,j}]^\top$, where $\begin{bmatrix}x_{r,j},\; y_{r,j}\end{bmatrix}$ is the radar's location.  
Given this augmented state, the radar's measurement model is
\begin{equation}
    h(\vect{x}_{l,i}, \vect{\theta}_{e,j}) = \begin{bmatrix}
    S_{E,i,j} \\ \phi_{i,j}
    \end{bmatrix} = \begin{bmatrix}
        \frac{P_{E,j}G_{I,i}\lambda^2}{(4\pi)^2 ||\vect{x}_{a,i}-\vect{x}_{r,j}||_2^2} \\ 
        \arctan\left(\frac{y_{r,j} - y_{l,i}}{x_{r,j}-x_{l,i}}\right)
    \end{bmatrix}.
\end{equation}

We assume that the measurements are corrupted with zero-mean Gaussian noise $\vect{\delta}$. If the position of the $i^\text{th}$ agent at the time of its $k^\text{th}$ measurement of the $j^\text{th}$ radar is $\vect{x}_{l,i}^k$ then the corresponding measurement is given by
\begin{equation}
\vect{z}_{i,j}^k = h(\vect{x}_{l,i}^k, \vect{\theta}_{e,j}) + \vect{\delta},\;
\vect{\delta}\sim \mathcal{N}(\vect{0},\Sigma_z),\; \Sigma_z = \begin{bmatrix}
    \sigma_{S_E}^2 & 0 \\
    0 & \sigma_{\phi}^2
\end{bmatrix},
\end{equation}
where $\sigma_{S_E}$ is the noise variance for the power measurement and $\sigma_\phi^2$ is the noise variance for the angle of arrival measurement. 

As agents traverse the environment, they gather measurements and store them in a set $\mylist{Z}_i = \{\vect{z}^k_{i,j}\}_{\forall k \in \{1,\hdots, N_{z,i}\}},$ where $N_{z,i}$ is the number of measurements agent $i$ has taken. The agents also store the locations where the measurements were taken, $\mylist{X}_{z,i} = \{\vect{x}_{i,j}^k\}_{\forall k \in \{1,\hdots,N_{z,i}\}}$.

We assume known data association, meaning the agents can reliably identify the radar station from which each measurement originates. This assumption allows us to focus on the core contribution of our work in path planning, rather than delving into the complexities of measurement data association. In real-world scenarios, agents can differentiate radar sources by leveraging signal characteristics such as frequency, pulse width, or pulse patterns, a topic explored in prior research~\cite{radarRecognition}.

We assume that the agents follow unicycle kinematics
\begin{equation} \label{eq:dynamics}
\begin{bmatrix}
\dot{x}(t) \\
\dot{y}(t) \\
\dot{\theta}(t)
\end{bmatrix} = 
\begin{bmatrix}
v(t)\cos{\theta(t)} \\
v(t)\sin{\theta(t)} \\
u(t)
\end{bmatrix},
\end{equation}
where $v(t)$ is the speed of the agent at time $t$ and $u(t)$ is the turn rate. Using the property of differential flatness, we can define kinematic feasibility constraints as done in~\cite{Buccieri2009_diffflat}. The velocity of the trajectory can be found as
\begin{equation}\label{eq:velocity_constraint}
    v(t) = ||\dot{\vect{p}}(t)||_2,
\end{equation}
and the turn rate $u(t)$ is 
\begin{equation}\label{eq:turn_rate_constraint}
u(t) = \frac{\dot{\vect{p}}(t) \times \ddot{\vect{p}}(t)}{{||\dot{\vect{p}}(t)||_2^2}},
\end{equation}
where $\vect{p}(t)=[x(t),y(t)]^\top$.
We can compute the curvature of the trajectory as
\begin{equation}\label{eq:curvature}
    \kappa(t) = \frac{u(t)}{v(t)}.
\end{equation}

 \section{Background} \label{sec:background}
 In this section, we describe several background topics that contribute to our algorithms. We first provide an overview of Voronoi diagrams, weighted Voronoi diagrams, and generalized Voronoi diagrams. Next, we discuss B-splines, which we use to parameterize paths.

\subsection{Voronoi Diagrams}
An \emph{ordinary Voronoi diagram} is defined using a set of generator points $\mylist{x}_g = \{\vect{x}_{g,1}, \hdots, \vect{x}_{g,N_g}\}$, where $N_g$ is the number of generator points and $\vect{x}_{g,i} \in \mathbb{R}^n$, with $n$ being the dimensions of the space (in our case, we use a planar space with $n=2$). 
The Voronoi cells are then defined as the region where the Euclidean distance between all points in the region and the generator point is less than the distance to any other generator point~\cite{boots1999spatial}:
\begin{equation}
    V(\vect{x}_{g,i}) = \{ \vect{x} \;|\; 
    \|\vect{x} - \vect{x}_{g,i} \|_2 \leq \|\vect{x} - \vect{x}_{g,j} \|_2,  
    \forall j \neq i ,\; j \in \{1,\hdots,N_g\} \}.
\end{equation}
The Voronoi diagram generated by the points $\mylist{x}_g$ is the set of all the Voronoi cells
\begin{equation}
    \mylist{V}(\mylist{x}_g) = \{V(\vect{x}_{g,1}),\hdots,V(\vect{x}_{g,N})\}.
\end{equation}
Voronoi edges are given as the intersection of two Voronoi regions if the intersection exists (e.g. 
 $V(\vect{x}_{g,i})\bigcap V(\vect{x}_{g,j})\neq \emptyset)$:
\begin{equation}
    e(\vect{x}_{g,i},\vect{x}_{g,j}) = \left(V(\vect{x}_{g,i})\bigcap V(\vect{x}_{g,j})\right).
\end{equation}
In the 2D case, these edges can either be line segments, half lines (where one direction of the line starts at a point and extends to infinity), or infinite lines. The endpoints of the Voronoi edges are called Voronoi vertices. These can also be defined as the intersection of three Voronoi regions~\cite{boots1999spatial}. We call the set of Voronoi edges $\mylist{E}(\mylist{x}_g)$ and the set of Voronoi vertices $\mylist{N}(\mylist{x}_g).$ 

Voronoi diagrams can be generated using various distance metrics. In a \textit{multiplicatively weighted Voronoi diagram}, the distance metric is modified by a weight associated with each generator point. The diagram is defined by the set of generator points $\mylist{x}_g$ and a corresponding set of positive weights $\mylist{w} = \{w_1, \dots, w_N\}$, where $w_i$ is the weight associated with the $i^\textit{th}$ generator. The multiplicatively weighted Voronoi region associated with generator $\vect{x}_{g,i}$ is defined as
\begin{equation}\label{eq:mult_voronoi}
    V(\vect{x}_{g,i}, w_i) = 
    \left\{ \vect{x} \;\middle|\; \frac{\|\vect{x} - \vect{x}_{g,i} \|_2}{w_i} 
    \leq \frac{\|\vect{x} - \vect{x}_{g,j} \|_2}{w_j},\ 
    \forall j  \neq i \right\}.
\end{equation}

The boundary between two such regions corresponds to a portion of an \textit{Apollonius circle}~\cite{boots1999spatial}—the set of points $\vect{x}$ that satisfy
\begin{equation}
\frac{\|\vect{x} - \vect{x}_{g,i}\|}{w_i} = \frac{\|\vect{x} - \vect{x}_{g,j}\|}{w_j}.
\end{equation}
Letting $\lambda = \frac{w_j}{w_i}$ denote the ratio of the weights, the Apollonius circle defined by points $\vect{x}_{g,i}$ and $\vect{x}_{g,j}$ has center $\vect{c} = (x_c, y_c)$ and radius $r$ given by

\begin{align}
    \vect{c} &= \frac{1}{1 - \lambda^2} \left( \vect{x}_{g,i} - \lambda^2 \vect{x}_{g,j} \right), \label{eq:apollonius_center} \\
    r &= \frac{\lambda}{|1 - \lambda^2|} \|\vect{x}_{g,i} - \vect{x}_{g,j}\|. \label{eq:apollonius_radius}
\end{align}
This circle represents the locus of points whose distances to $\vect{x}_i$ and $\vect{x}_j$ are in a fixed ratio equal to the inverse ratio of the weights. The corresponding Voronoi edge is a portion of this circle that lies within both weighted Voronoi regions. When $w_i = w_j$, the Apollonius circle degenerates to a straight line, which is the classical Voronoi boundary.

An edge between two Voronoi regions exists where their weighted regions intersect:
\begin{equation}
    e(\vect{x}_{g,i}, \vect{x}_{g,j}) = V(\vect{x}_{g,i}, w_i) \cap V(\vect{x}_{g,j}, w_j).
\end{equation}

The \textit{Voronoi vertices} of the diagram are the points where three or more such edges intersect. We denote the set of edges as $\mylist{E}(\mylist{x}_g, \mylist{w})$ and the set of vertices as $\mylist{N}(\mylist{x}_g, \mylist{w})$. The full set of edges and vertices in a multiplicatively weighted Voronoi diagram can be computed efficiently using the algorithm described in~\cite{held2020efficient}.

\subsection{B-splines}
B-splines are piecewise polynomial functions defined by control points $\mylist{c} = \{\vect{c}_1,\hdots,\vect{c}_{N_c}\}, \vect{c}_i\in\mathbb{R}^2$ and knot points $\vect{t}_k = \{t_0-p\Delta_t,\hdots,t_0-\Delta_t, t_0,t_0+\Delta_k,t_0+2\Delta_k,\hdots t_f,t_f+\Delta_k,\hdots t_f+p\Delta_k\}$ where $\Delta_k = (t_f-t_0)/(N_k-2p)$ is the knot point spacing, $N_k = N_c+p+1$ is the number of knot points and $p$ is the degree of the B-spline for an unclamped uniform B-spline. The B-spline is defined on the interval $[t_0,t_f]$ as
\begin{equation} \label{eq:b-spline_basis}
\vect{p}(t) = \sum ^{{N_c}}_{i=1} B_{i,p}(t)\vect{c}_i,
\end{equation}
where the basis functions $B_{i,p}$ are defined using the Cox-de Boor recursive formula shown in~\cite{cox1972numerical}. B-splines are commonly used in path planning applications because of their local support property (sparse Jacobians) and the convex hull property (the trajectory must be within the convex hull of the control points)~\cite{judd2001spline}.

\section{Radar Localization and Probability of Detection}\label{sec:radar_localization_PD}
In this section, we outline our method for localizing radar and estimating the ERP using angle-of-arrival and signal strength measurements. Our approach to radar localization utilizes a combination of nonlinear least squares and an EKF to track enemy radar. We also present a method similar to~\cite{costley2022sensitivity,costley2023sensitivity} that accounts for uncertainty in the radar parameters to find the radar PD. This approach linearizes the radar PD equations and propagates the uncertainty through the linearized model. 
\subsection{Radar Localization} \label{sec:radar_localization}
Our approach for localizing radar stations uses a non-linear least squares method to initialize the radar station models, leveraging available measurements to estimate their initial parameters. 
Once a model is initialized, the state estimates are refined through EKF correction updates as new measurements are received. 
The goal is to estimate each radar's augmented state $\vect{\theta}_{e,j}$; its location and ERP, and the associated uncertainty of that estimate, represented as a covariance matrix. 
The estimation algorithm outputs a list of mean values and covariance pairs, 
$\mylist{R} = \{(\mu_{\vect{\theta}_{e,j}}, \Sigma_{\vect{\theta}_{e,j}})\}\;\forall j \in \{1,\hdots,N_{\hat{r}}\},$where $N_{\hat{r}}$ is the total number of estimated radar stations. 
Algorithm~\ref{alg:EKF_RANSAC} provides a detailed outline of the radar estimation process.

The input to the algorithm is a stream of measurements $\vect{z}_{i,j}^k$ and measurement locations $\vect{x}_{i,j}^k$ from all agents. The Algorithm outputs a mean and covariance estimate for each radar station that has been discovered $\mylist{R}$. On line~\ref{alg:EKF_initialize}, several lists are initialized; $\mylist{R}$ will store the mean and covariance pairs for each radar station, then for each radar station a list to store the measurements $\mylist{Z}_j$ and measurement locations $\mylist{x}_j$ are created. 
When a new measurement $\vect{z}^k_{i,j}$ and measurement location $\vect{x}^k_{i,j}$ is received, it is added to the list of measurements $\mylist{Z}_j$ and measurement locations $\mylist{x}_j$ for each radar on line~\ref{alg:ekf_add_measurement}. 

If a model (mean and covariance estimate) exists for the $j^\textit{th}$ radar, 
the algorithm uses the EKF measurement correction step to incorporate the measurement into the existing model (lines~\ref{alg:ekf_check_model}-\ref{alg:ekf_mean_update}). This is done by computing the Kalman gain using the current estimate covariance for the $j^{\text{th}}$ radar $\Sigma_{\vect{\theta}_{e,j}},$ the measurement covariance $\Sigma_z$, and the Jacobian of the measurement model 
\begin{equation}
J_h(\vect{x}_{i,j}^k, \mu_{\vect{\theta}_{e,j}}) = \left. \frac{\partial h }{ \partial \vect{\theta}_{e,j}} \right|_{(\vect{x}_{l,i}^k, \vect{\theta}_{e,j}) = (\vect{x}_{i,j}^k, \mu_{\vect{\theta}_{e,j}})}
\end{equation}
evaluated at the measurement location $\vect{x}_{i,j}^k$ and the current mean value $\mu_{\vect{\theta}_{e,j}}$.  After the Kalman gain is computed, the mean and covariance values are updated.

If a model does not exist for a specific radar, the algorithm will initialize a new one if there are $N_{z,\text{min}}$ measurements (lines~\ref{alg:ekf_check_num}-\ref{alg:ekf_add_model}). 
A new track's mean value $\mu_{\vect{\theta}_{e,j}}$ is calculated using a non-linear least squares optimization algorithm where the objective function is to minimize the sum of the squared Mahalonobis distances between the measurements and the measurement model. 
The covariance is approximated by finding the inverse of the Fisher information matrix. 
To do this, we find the Jacobian of each measurement model with respect to the measurement location $\vect{x}_{i,j}^k$ and stack them as
\begin{equation}\label{eq:Jac_stacked}
    \boldsymbol{J}_h = 
    \begin{bmatrix}
       J_h(\vect{x}_{i,j}^1, \mu_{\vect{\theta}_{e,j}})\\
       \vdots\\
       J_h(\vect{x}_{i,j}^k, \mu_{\vect{\theta}_{e,j}})\\
       \vdots \\
       J_h(\vect{x}_{i,j}^{N_{z,\textit{min}}}, \mu_{\vect{\theta}_{e,j}})
    \end{bmatrix}\forall\vect{x}_{i,j}^k\in\mylist{x}_j,\in\mathbb{R}^{3\times2N_{z,min}}
\end{equation}
where $i$ represents the agent index of the agent that took each measurement (potentially different).
Similarly, the measurement covariance is stacked in a block diagonal form as 
\begin{equation}\label{eq:cov_stacked}
\boldsymbol{\Sigma}_z = 
\begin{bmatrix}
\Sigma_z & \mathbf{0} & \mathbf{0} & \cdots & 0 \\
\mathbf{0} & \Sigma_z & \mathbf{0} & \cdots & 0 \\
\mathbf{0} & \mathbf{0} & \Sigma_z & \cdots & 0 \\
\vdots & \vdots & \vdots & \ddots & \vdots \\
0 & 0 & 0 & \cdots & \Sigma_z
\end{bmatrix}\in \mathbb{R}^{2N_{z,min}\times2N_{z,min}}.
\end{equation}
Equations~\eqref{eq:Jac_stacked} and~\eqref{eq:cov_stacked} are then used to approximate the covariance of the estimate of the radar's location and ERP (line~\ref{alg:ekf_FIM}). The new model $(\mu_{\vect{\theta}_{e,j}},\Sigma_{\vect{\theta}_{e,j}})$ is added to the list of models $\mylist{R}$ in line~\ref{alg:ekf_add_model}.

\begin{algorithm} 
 \caption{EKF for radar model estimation}\label{alg:EKF_RANSAC}
 \begin{algorithmic}[1]
 \STATE{\textbf{Input:} stream of measurement $\vect{z}_{i,j}^k$ and measurement locations $\vect{x}_{i,j}^k$ from all agents}
 \STATE{\textbf{Output:} $\mylist{R}$}
\STATE{\textbf{Initialize:} $\mylist{R}\leftarrow\emptyset,\mylist{Z}_j\leftarrow\emptyset,\mylist{X}_j\leftarrow\emptyset, N_{z,j} = 0 \;\forall j \in \{1,\hdots,N_r\}$,Initialized$_j$ $\leftarrow$ False  $\quad \forall j \in \{1, \dots, N_r\}$
}\label{alg:EKF_initialize}
\FOR{$\vect{z}_{i,j}^k,\vect{x}^k_{i,j}$ in stream}\label{alg:start_EKF}
\STATE{$\mylist{x}_j=\mylist{x}_j\bigcup\{\vect{x}_{i,j}^k\},\mylist{z}_j=\mylist{z}_j\bigcup \{\vect{z}_{i,j}^k\}$}\label{alg:ekf_add_measurement}
\IF{Initialized$_j$}\label{alg:ekf_check_model}
\STATE{$K = \Sigma_{\vect{\theta}_{e,j}} J_h(\vect{x}_{i,j}^k,\mu_{\vect{\theta}_{e,j}})^\top$}
\STATE{\hspace{1.5em}$\cdot \left(J_h(\vect{x}_{i,j}^k,\mu_{\vect{\theta}_{e,j}}) 
       \Sigma_{\vect{\theta}_{e,j}} 
       J_h(\vect{x}_{i,j}^k,\mu_{\vect{\theta}_{e,j}})^\top 
       + \Sigma_z \right)^{-1}$}\label{alg:Kalman_gain}
\STATE{$\mu_{\vect{\theta}_{e,j}} = \mu_{\vect{\theta}_{e,j}} + K(\vect{z}_i^k-h(\vect{x}_{i,j}^k,\mu_{\vect{\theta}_{e,j}}))$}\label{alg:ekf_cov_update}
\STATE{$\Sigma_{\vect{\theta}_{e,j}} = (I-KJ_h)\Sigma_{\vect{\theta}_{e,j}}$}\label{alg:ekf_mean_update}
\ELSIF{$|\mylist{z}_j| = N_{r,\textit{min}}$}\label{alg:ekf_check_num}
\STATE{$\mu_{\vect{\theta}_{e,j}} = \operatornamewithlimits{argmin}\limits_{\mu_{\vect{\theta}_{e,j}}} \sum\limits_{\vect{x}_{i,j}^k \in \mylist{X}_j, \vect{z}_j^k \in \mylist{Z}_j} ||h(\vect{x}_{i,j}^k, \mu_{\vect{\theta}_{e,j}}) \!-\! \vect{z}_j^k||^2_{\Sigma_z}$}\label{alg:nonlinear_least_squares}
\STATE{$\Sigma_{\vect{\theta}_{e,j}} = (\boldsymbol{J}_h\boldsymbol{\Sigma}_z^{-1}\boldsymbol{J}_h^\top)^{-1}
$}\label{alg:ekf_FIM}
\STATE{$\mylist{R} \leftarrow \mylist{R} \bigcup (\mu_{\vect{\theta}_{e,j}},\Sigma_{\vect{\theta}_{e,j}})$}\label{alg:ekf_add_model}
\STATE{Initialized$_j\leftarrow$ True}
\ELSE
\STATE{Pass}
\ENDIF
\ENDFOR
 \end{algorithmic} 
 \end{algorithm}

\subsection{Radar Probability of Detection} \label{sec:RADAR_PD}
We now present the method used to calculate PD and its associated uncertainty, employing a linearization technique similar to the approach described in~\cite{costley2023sensitivity}.
The SNR and PD are functions of known parameters of the agent $\vect{\theta}_{k,i} = [\sigma_i,\vect{x}_i]$ (radar cross section and position of the agent), unknown parameters of the radar $\vect{\theta}_{u,j} = [P_{\textit{fa},j}, G_{R,j}, \lambda, \tau_{p,j} , T_{s,j}]^\top$
(probability of false alarm, radar receive gain, wavelength, radar pulse length, radar system temperature)  and estimated parameters of the radar $\vect{\theta}_{e,j}$ (radar position and ERP). The known parameters can be found before the mission starts (radar cross section), or found by the agent during the mission (we assumed known locations). The unknown parameters are unobservable by the agents using the measurements described in Section~\ref{sec:Problem_statement}. Instead, we assume the agents have some knowledge (a probability distribution) of these parameters. The estimated parameters are found using the method described in Section~\ref{sec:radar_localization}. 

The signal-to-noise ratio (SNR) from $i^\textit{th}$ agent to the $j^\textit{th}$ radar is
\begin{equation}
        \mathrm{SNR}_{i,j}(\vect{\theta}_{k,i}, \vect{\theta}_{u,j}, \vect{\theta}_{e,j}) = \frac{P_{E,j}G_{R,j}\lambda^2\sigma_i \tau_{p,j}}{(4\pi)^3R_{i,j}^4\kappa T_{s,j}L_j},
\end{equation}
where $G_{R,j}$ is the gain of the receive antennae of the radar, $\sigma_i$ is the radar cross section (RCS) of the agent, $\tau_{p,j}$ is the radar pulse length, $\kappa$ is the Boltzman constant, and $T_{s,j}$ is the radar system noise. Uing the SNR the PD of the $i^\textit{th}$ agent from the $j^\textit{th}$ radar is
\begin{equation}
    P_{D,i,j}(\vect{\theta}_{k,i}, \vect{\theta}_{u,j}, \vect{\theta}_{e,j}) {=} \exp{\left(\frac{\ln{P_{\textit{fa},j}}}{\mathrm{SNR}_{i,j}(\vect{\theta}_{k,i}, \vect{\theta}_{u,j}, \vect{\theta}_{e,j}){+}1}\right)},
\end{equation}
where $P_{\textit{fa},j}$ is the probability of a false alarm (determined by the radar operator) and $\mathrm{SNR}_{i,j}$ is the signal-to-noise ratio. 

Using the current estimate of the radar parameters, the overall PD of the $i^\textit{th}$ agent is given as
\begin{equation} \label{eq:overall_PD}
    P_{D,i}(\vect{\theta}_{k,i},\mylist{\Theta}_{u},\mylist{\Theta}_e) = 1{-}\prod_{j=1}^{N_{\hat{r}}}\left(1{-}P_{D,i,j}(\vect{\theta}_{k,i}, \vect{\theta}_{u,j}, \vect{\theta}_{e,j})\right),
\end{equation}
where $N_{\hat{r}}$ is the current number of discovered radar, $\mylist{\Theta}_u = \{\vect{\theta}_{u,j}\}_{\forall j \in \{1,\hdots,N_r\}}$ is the set of unknown radar parameters for radar systems within range of the $i^\textit{th}$ agent, and $\mylist{\Theta}_e = \{\vect{\theta}_{e,j}\}_{\forall j \in \{1,\hdots,N_r\}}$ is the set of estimated radar parameters also for radar systems within range of the $i^\textit{th}$ agent.

We assume that we have a known Gaussian distribution $\vect{\theta}_{k,i} \sim \mathcal{N}(\mu_{\vect{\theta}_{k,i}},\Sigma_{\vect{\theta}_{k,i}})$ that quantifies the uncertainty in the agent's known parameters. To get a similar distribution for the radar's estimated parameters, we use the estimator described in Section~\ref{sec:radar_localization} to obtain mean value and covariance estimates for each discovered radar $\mylist{R}$. We then combine these estimated parameters into a mean vector
\begin{equation}
    \mu_{\mylist{\Theta}_{e}} = 
    \begin{bmatrix}
     \mu_{\vect{\theta}_{e,1}}\\  
     \mu_{\vect{\theta}_{e,2}}\\  
     \vdots\\
     \mu_{\vect{\theta}_{e,N_{\hat{r}}}}\\  
    \end{bmatrix}
\end{equation}
and covariance matrix
\begin{equation}
    \Sigma_{\mylist{\Theta}_{e}} = 
    \begin{bmatrix}
     \Sigma_{\vect{\theta}_{e,1}} & \boldsymbol{0}&  \hdots &    \boldsymbol{0} \\
     \boldsymbol{0} &\Sigma_{\vect{\theta}_{e,2}} &  \hdots &    \boldsymbol{0} \\
     \vdots & \hdots  & \ddots &    \vdots \\
      \boldsymbol{0} & \hdots & \hdots & \Sigma_{\vect{\theta}_{e,N_{\hat{r}}}}  \\
    \end{bmatrix}\in \mathbb{R}^{5N_{\hat{r}}\times5N_{\hat{r}}},
\end{equation}
where the individual radar estimates are uncorrelated and the combined radar parameters are normally distributed $\mylist{\Theta}_{e} \sim \mathcal{N}(\mu_{\mylist{\Theta}_{e}},\Sigma_{\mylist{\Theta}_{e}})$. 

Because we cannot estimate the unknown parameters, we must use prior knowledge of the system to create reasonable prior beliefs of these parameters. We assume a normal distribution and represent the mean value of this prior belief as 
\begin{equation}
    \mu_{\mylist{\Theta}_{u}} = 
    \begin{bmatrix}
     \mu_{\vect{\theta}_{u,1}}\\  
     \mu_{\vect{\theta}_{u,2}}\\  
     \vdots\\
     \mu_{\vect{\theta}_{u,N_{\hat{r}}}}\\  
    \end{bmatrix}
\end{equation}
and the covariance matrix as
\begin{equation}
    \Sigma_{\mylist{\Theta}_{u}} = 
    \begin{bmatrix}
     \Sigma_{\vect{\theta}_{u,1}} & \boldsymbol{0}&  \hdots &    \boldsymbol{0} \\
     \boldsymbol{0} &\Sigma_{\vect{\theta}_{u,2}} &  \hdots &    \boldsymbol{0} \\
     \vdots & \hdots  & \ddots &    \vdots \\
      \boldsymbol{0} & \hdots & \hdots & \Sigma_{\vect{\Theta}_{u,N_{\hat{r}}}}  \\
    \end{bmatrix}\in \mathbb{R}^{3N_{\hat{r}}\times3N_{\hat{r}}},
\end{equation}
where $\mu_{\vect{\theta}_{u,j}}$ and $\Sigma_{\vect{\theta}_{u,j}}$ are the mean value and covariance for the unknown parameters of the $j^{\textit{th}}$ radar station and the combined unknown parameters are normally distributed $\mylist{\Theta}_{u} \sim \mathcal{N}(\mu_{\mylist{\Theta}_{u}},\Sigma_{\mylist{\Theta}_{u}})$. In this work, we use the same prior distribution for each radar's unknown parameters ($\mu_{\vect{\theta}_{u,j}}=\mu_{\vect{\theta}_{u,i}},\Sigma_{\vect{\theta}_{u,j}}=\Sigma_{\vect{\theta}_{u,i}},\forall (i,j) \in \{1,\hdots,N_{\hat{r}}\}).$

Using the known, unknown, and estimated parameter distributions, we wish to compute the PD and provide a covariance for that estimate. We do this by creating a first-order approximation of the PD (Equation~\eqref{eq:overall_PD}) centered about the means,
\begin{equation} \label{eq:first_order_approx}
P_{D,i}(\mu_{\vect{\theta}_{k,i}}+\vect{\delta}_k,\mu_{\mylist{\Theta}_{u}} + \vect{\delta}_u,\mu_{\mylist{\Theta}_e}+\vect{\delta}_e) \approx \\
P_{D,i}(\mu_{\vect{\theta}_{k,i}},\mu_{\mylist{\Theta}_{u}},\mu_{\mylist{\Theta}_e}) + \vect{\delta}_k J_k+  
 \vect{\delta}_u J_u 
+ \vect{\delta}_e J_e,
\end{equation}
where $\vect{\delta}_k\in\mathbb{R}^3, \vect{\delta}_u \in \mathbb{R}^{5N_{\hat{r}}}$, and $\vect{\delta}_e \in \mathbb{R}^{3N_{\hat{r}}},$  are perturbations in known, unknown, and estimated parameters respectively.  The Jacobian of Equation~\eqref{eq:overall_PD} with respect to the agent's parameters is $J_k={\partial P_{D,i}(\mu_{\vect{\theta}_{k,i}},\mu_{\mylist{\Theta}_{u}},\mu_{\mylist{\Theta}_e})}/{\partial \vect{\theta}_{k,i}}\in \mathbb{R}^{1\times3}$.  And similarly, $J_u={\partial P_{D,i}(\mu_{\vect{\theta}_{k,i}},\mu_{\mylist{\Theta}_{u}},\mu_{\mylist{\Theta}_e})}/{\partial \mylist{\Theta}_{u}}\in \mathbb{R}^{1\times5N_{\hat{r}}}$ is the Jacobian of Equation~\eqref{eq:overall_PD} with respect to the unknown radar parameters $\mylist{\Theta}_{u},$ and $J_e{\partial P_{D,i}(\mu_{\vect{\theta}_{k,i}},\mu_{\mylist{\Theta}_{u}},\mu_{\mylist{\Theta}_e})}/{\partial \mylist{\Theta}_{e}}\in \mathbb{R}^{1\times3N_{\hat{r}}}$ is the Jacobian of Equation~\eqref{eq:overall_PD} with respect to the estimated radar parameters $\mylist{\Theta}_{e}.$ These Jacobians can be found analytically or through automatic differentiation. In this work, we use the JAX automatic differentiation library~\cite{jax2018github}.

Using these distributions and the linearized model of PD we can find the mean value of the PD as
\begin{equation} \label{eq:PD_mean}
    \mu_{P_{D,i}}(\vect{\theta}_{k,i},\mylist{\Theta}_u,\mylist{\Theta}_e) = P_{D,i}(\mu_{\vect{\theta}_{k,i}},\mu_{\mylist{\Theta}_{u}},\mu_{\mylist{\Theta}_{e}}).
\end{equation}
and assuming the known, unknown, and estimated parameters are independent, the variance
\begin{equation}\label{eq:PD_cov}
    \sigma_{P_{D,i}}^2(\vect{\theta}_{k,i},\mylist{\Theta}_u,\mylist{\Theta}_e)    =\\  J_k\Sigma_{\vect{\theta}_{k,i}}J_k^\top+  
    J_u\Sigma_{\mylist{\Theta}_{u}}J_u^\top+  
    J_e\Sigma_{\mylist{\Theta}_{e}}J_e^\top.
\end{equation}
Using this mean and covariance, we can create an approximate distribution for the PD, $P_{D,i} \sim \mathcal{N}(\mu_{P_{D,i}},\sigma^2_{P_{D,i}}).$ We can then use this distribution to approximate the probability of the true PD being below a threshold:
\begin{equation}\label{eq:prob_pd_less_than_threshold}
    P(P_D \leq P_{D,t}) =\\ 
    \Phi\left(P_{D,t},\mu_{P_D}(\vect{\theta}_{k,i},\mylist{\Theta}_u,\mylist{\Theta}_e),\sigma^2_{P_D}(\vect{\theta}_{k,i},\mylist{\Theta}_u,\mylist{\Theta}_e)\right),
\end{equation}
where 
\begin{equation}
\Phi(x; \mu, \sigma^2) = \frac{1}{2} \left[ 1 + \mathrm{erf}\left( \frac{x - \mu}{\sqrt{2} \sigma} \right) \right] 
\end{equation}
is the normal cumulative distribution function and
\begin{equation}
    \mathrm{erf}(z) = \frac{2}{\sqrt{\pi}} \int_0^z e^{-t^2} \, dt.
\end{equation}

\section{Low-Priority Agent Path Planning}\label{sec:lp-path-planning}
The objective of low-priority agents is to explore the operational area, discover unknown enemy radar, and find a safe path for the high-priority agent. We employ a decentralized optimization strategy. Each agent uses the current information it has received about the environment to plan its own path.  This path information is transmitted to other agents, who then plan their paths. The paths are planned in a round-robin fashion, where the agent closest to the most uncertain radar station plans their path first. 

\subsection{Objective Function} \label{sec:LPP_objective}
The agents' objectives of discovering enemy radar, improving radar parameter estimates, and finding the optimal trajectory for the high-priority agent lead to a three-part objective function. First, agents are incentivized to \emph{explore} previously unexplored regions through a Bayesian approach, detailed later in Equation~\eqref{eq:exploration_objective}. Second, agents aim to \emph{reduce the uncertainty} in radar parameter estimates by selecting waypoints that improve the Extended Kalman Filter (EKF) estimation accuracy, as quantified in Equation~\eqref{eq:uncertainty_objective}. Third, to efficiently reach its goal, the high-priority agent seeks a trajectory close to the \emph{straight-line path} between its start location, $\vect{x}_{h_0}$, and goal location, $\vect{x}_{h_f}$. Deviations from this path are penalized, encouraging exploration near the optimal straight-line route, as shown in Equation~\eqref{eq:strait_line_distance_objective}.

Waypoints for the low-priority agents are generated using IPOPT~\cite{IPOPT}, a nonlinear interior point optimization algorithm. The optimization balances the trade-off between exploiting existing knowledge and exploring new areas via the three aforementioned objectives:
\begin{enumerate}
    \item \emph{Exploration},
    \item \emph{Radar parameter uncertainty reduction}, and
    \item \emph{Goal-directed prioritization of the high-priority agent's trajectory}
\end{enumerate}
These objectives and their implementation in determining low-priority vehicle paths are described in detail below.
\begin{figure*}[ht]
    \centering
\includegraphics[width=0.95\linewidth,trim={.62cm .2cm 0.2cm .1cm},clip]{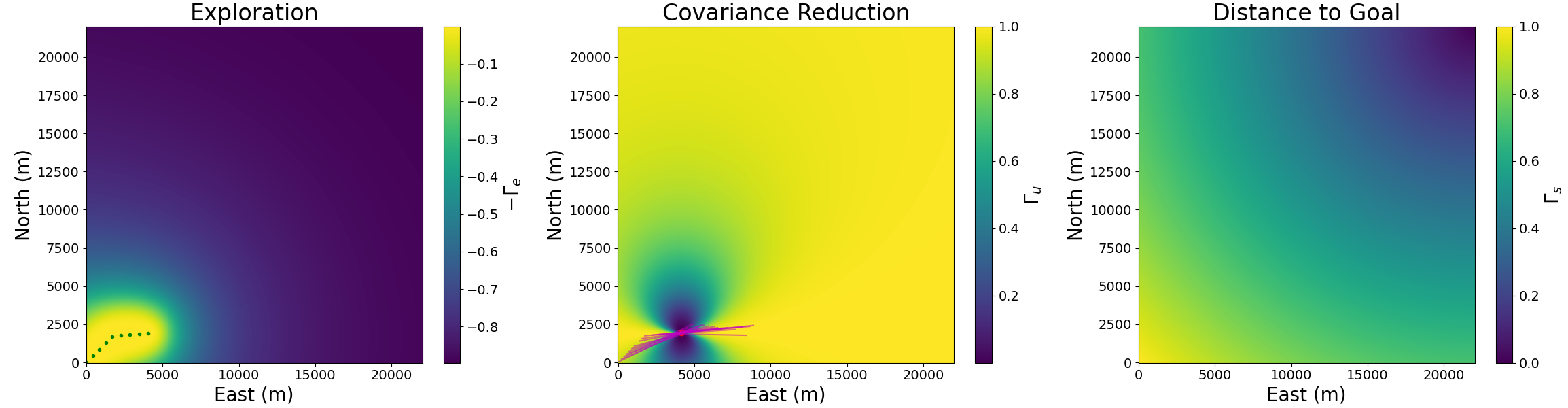}
    \caption{This figure shows the three different parts of the objective function. The left figure illustrates the exploration objective function $\Gamma_e(\vect{x})$, where the green points show locations the agents have already explored. We show the negative of $\Gamma_e(\vect{x})$ because it is a penalty, so the darker color represents better measurement locations. Going further from the explored area provides better measurements. The middle figure shows the covariance reduction objective function $\Gamma_u(\vect{x})$. The estimated radar location is plotted in red, and the previous measurements are shown as magenta lines. From this figure, we see that the best measurement locations to reduce a radar's covariance uncertainty, according to this objective function, come from being closer to the radar and from going perpendicular to the current measurements.  The final plot shows the distance to the goal, where the goal is in the upper right corner of the plot.}\label{fig:lpp_objectives}
\end{figure*}

\subsubsection{Exploration}
Agents are encouraged to visit locations that have not yet been explored.  
To do this, each agent maintains a list of previously explored locations:
\begin{equation}
\mylist{X}_{e,i} = \\\{\vect{x}_{l,i}(t_0),\,\vect{x}_{l,i}(t_0+\Delta_{t_e}),\dots,\,\vect{x}_{l,i}(t_0+N_e\Delta_{t_e})\},
\end{equation}
where $\vect{x}_{l,i}(t)$ is the state of the $i^\text{th}$ low-priority agent at time $t$, $t_0$ is the mission start time, $\Delta_{t_e}$ is the time interval between stored path points, and $N_e$ is the number of points stored up to the current time $t_c$, defined as:
\begin{equation}
N_e = \text{int}\left(\frac{t_c - t_0}{\Delta_{t_e}}\right).
\end{equation}
The individual path history list of every agent is combined into a single list $\mylist{X}_e = \bigcup_{i=1}^{N_l}\mylist{x}_{e,i}$.

From this combined path history, we will approximate the probability that an undiscovered radar is present at a location $\vect{x}$ given the locations that the agents have visited.
We do this using Bayes' rule and approximating the probability that an agent intercepts an enemy radar signal.
Given an agent was at location $\vect{x}_{l,i}(t_j)$, and that a radar is at location $\vect{x}$, the probability the agent would have intercepted a signal is given by:
\begin{equation}
    P(\textit{intercept at }\vect{x}_{l,i}(t_j)\mid\textit{radar at }\vect{x}) 
    = \exp\left(\frac{\ln(P_{fa})}{\mathrm{SNR}\left(\vect{x}_{l,i}(t_j),\vect{x}\right)+1}\right),
\end{equation}
where $P_{fa}$ is the probability of false alarm, and $\mathrm{SNR}(\vect{x}_{i},\vect{x})$ is the signal-to-noise ratio (SNR) between the agent and the radar. The SNR is given by:
\begin{equation}
    \mathrm{SNR}(\vect{x}_{l,i}(t_j),\vect{x}) = \frac{P_{T}G_{T}G_{I,i}\lambda^2\tau_{p}}{(4\pi)^2||\vect{x}_{l,i}(t_j)-\vect{x}||_2^2\kappa T_{s}L\delta_i}.
\end{equation}
Since the parameters of the potentially undiscovered radar are unknown, we assume default values for the transmitted power, $P_T$, the transmit gain, $G_T$, the pulse width, $\tau_p$, the system temperature, $T_s$, and the loss factor, $L$. 
However, we have prior knowledge of the agent's intercept antenna gain, $G_{I,i}$, and the signal wavelength, $\lambda$. 
In addition, we introduce a discount term, $\delta_l$, to account for signal degradation caused by factors such as misaligned antennae, mismatched filters, and other unknown losses. 
The probability of failing to intercept a signal at $\vect{x}_i$, given that a radar is located at $\vect{x}$, is related to the probability of intercepting by $P(\textit{no intercept at } \vect{x}_i | \textit{radar at } \vect{x}) = 1 - P(\textit{intercept at } \vect{x}_i | \textit{radar at } \vect{x}).$

We wish to approximate the probability that a radar could be found at location $\vect{x}$ given a history that the agents did not intercept signals at the locations $\mylist{X}_e$.  Using Bayes' rule, this is given by:
\begin{equation}\label{eq:bayes}
    P\left(\textit{radar at } \vect{x} \,\middle|\, \bigcap_{\vect{x}_i \in \mylist{X}_e} \textit{no intercept at } \vect{x}_i\right) = \\
    \frac{P\left(\bigcap_{\vect{x}_i \in \mylist{X}_e} \textit{no intercept at } \vect{x}_i \,\middle|\, \textit{radar at } \vect{x}\right) P(\textit{radar at } \vect{x})}
    {P\left(\bigcap_{\vect{x}_i \in \mylist{X}_e} \textit{no intercept at } \vect{x}_i\right)}.
\end{equation}
If we assume that intercepting a signal at $\vect{x}_i$ is independent of intercepting a signal at a different location, given that there is a radar at $\vect{x}$, then  
\begin{equation}\label{eq:cond_prob}
    P\left(\bigcap_{\vect{x}_i \in \mylist{X}_e} \textit{no intercept at } \vect{x}_i \,\middle|\, \textit{radar at } \vect{x}\right)
    = \prod_{\vect{x}_i \in \mylist{X}_e} P\left(\textit{no intercept at } \vect{x}_i \,\middle|\, \textit{radar at } \vect{x}\right).
\end{equation}

We assume a uniform prior radar distribution $\Phi(\vect{x})=P(\textit{radar at }\vect{x}) = 0.5$, meaning there are equal chances that a radar is present at $\vect{x}$ or not. If additional information about radar locations in the region were available, it could be used to refine this prior probability. To find the denominator, we use a partition:  
\begin{multline}\label{eq:total_prob}
    P\left(\bigcap_{\vect{x}_i \in \mylist{X}_e} \textit{no intercept at } \vect{x}_i\right) = 
    P\left(\bigcap_{\vect{x}_i \in \mylist{X}_e} \textit{no intercept at } \vect{x}_i \,\middle|\, \textit{radar at } \vect{x}\right) P(\textit{radar at } \vect{x}) + \\
    P\left(\bigcap_{\vect{x}_i \in \mylist{X}_e} \textit{no intercept at } \vect{x}_i \,\middle|\, \textit{no radar at } \vect{x}\right) P(\textit{no radar at } \vect{x}).
\end{multline}

We know $P(\textit{no intercept at }\vect{x}_i | \textit{no radar at }\vect{x}) = 1 - P_{f_a}$. Assuming that not intercepting a signal at one location is independent of not intercepting a signal at another location when no radar is present, we can write  
$P(\bigcap_{\vect{x}_i \in \mylist{X}_e} \textit{no intercept at }\vect{x}_i | \textit{no radar at }\vect{x}) = (1 - P_{f_a})^{|\mylist{x}_e|}$, where $|\mylist{x}_e|$ is the cardinality of $\mylist{x}_e$.

From this, we compute the probability that an undiscovered radar exists at a location $\vect{x}$ by substituting Equations~\eqref{eq:total_prob} and~\eqref{eq:cond_prob} into Equation~\eqref{eq:bayes}
\begin{equation} \label{eq:exploration_objective}
\Gamma_e(\vect{x},\mylist{x}e) = 
\frac{\prod\limits_{\vect{x}_i \in \mylist{x}e}\left(1 - \exp\left(\frac{\ln(P_{f_a})}{\mathrm{SNR}(\vect{x}_{i},\vect{x}) + 1}\right)\right) \Phi(\vect{x})}{\prod\limits_{\vect{x}_i \in \mylist{x}e}\left(1 - \exp\left(\frac{\ln(P_{fa})}{\mathrm{SNR}(\vect{x}_{i},\vect{x}) + 1}\right)\right)\Phi(\vect{x}) + (1 - P_{f_a})^{|\mylist{x}_e|}(1 - \Phi(\vect{x}))}.
\end{equation}
We aim to maximize the probability of discovering previously undetected radars by directing low-priority agents to locations where $\Gamma_e(\vect{x},\mylist{x}_e)$ is high. The effect of this objective term is shown in the left panel of Figure~\ref{fig:lpp_objectives}, where the agent's path history is visualized in green. As shown, this component of the objective encourages agents to move toward areas that are far from previously visited locations.

\subsubsection{Radar Parameter Uncertainty Reduction}
The second part of the objective function, $\Gamma_u(\vect{x})$, guides agents toward areas where radar measurements are more informative and will improve localization. We do this by noting that when a measurement is incorporated into an EKF, the covariance does not depend on the value that is measured, only on the location at which the measurement was taken. 
This component of the objective function incentivizes agents to take measurements at locations that reduce the determinant of the radar parameter covariance matrix, similar to the approach in~\cite{xu2024systematical}, where the authors minimize the determinant of an EKF covariance matrix to improve parameter estimation accuracy using angle-of-arrival measurements.

From the EKF described in Section~\ref{sec:radar_localization}, we have a list of the current estimated mean and covariance values for the radar parameters $\mylist{R}$. The equation for the updated covariance when incorporating an additional measurement is
\begin{equation}
    \operatorname{cov}^+(\vect{x},\mu_{\vect{x}_{r,j}},\Sigma_{\vect{x}_{r,j}})= (I-KJ_h(\vect{x},\mu_{\vect{x}_{r,j}}))\Sigma_{\bar{\vect{x}}_{r,j}}, 
\end{equation}
where $K$ is the Kalman gain defined in Algorithm~\ref{alg:Kalman_gain}, $I$ is the identity matrix, and $J_h$ is the Jacobian of the measurement model. Given a list of models--each defined by their mean and covariance--for all the discovered radars, we aim to find the measurement location that would most improve all the models.  To do this, we compute the mean of the determinant of all the updated covariances that would result if a measurement were taken at location $\vect{x}$. This yields, 
\begin{equation}\label{eq:uncertainty_objective}
    \Gamma_u(\vect{x},\mylist{R}) =\\ \frac{1}{N_{\hat{r}}}\sum\limits_{(\mu_{\vect{x}_{r,j}},\Sigma_{\vect{x}_{r,j}})\in\mylist{R}}\frac{\operatorname{det}(\operatorname{cov}^+(\vect{x},\mu_{\vect{x}_{r,j}},\Sigma_{\vect{x}_{r,j}}))}{d_\textit{cov}},
\end{equation}
where $d_\textit{cov}$ is a normalizing term that controls how much the objective function prioritizes reducing the covariance. Once the mean value of the determinant of the future covariances is much less than $d_\textit{cov}$, the covariance reduction objective will no longer affect the objective function. This allows agents to ignore radar stations that are localized well enough, where well enough is defined by $d_\textit{cov}$.
At location $\vect{x}$, it is not likely that a measurement will be received for each model; however, this objective is a good heuristic to see how traveling to a certain location will improve the current estimates of the enemy radar.
This part of the objective function is visualized in the center panel of Figure~\ref{fig:lpp_objectives}. Past angle of arrival measurements are shown in magenta, and the heatmap represents the corresponding objective function values. As illustrated, taking new measurements at locations that are off-axis from the previous ones results in the greatest reduction in uncertainty.

\subsubsection{Goal-directed Search Prioritization}
The final piece of the objective function, $\Gamma_s(\vect{x})$, incentivizes the low-priority agents to explore in the direction of the high-priority agent's goal. The objective is the normalized distance to the high-priority goal, 
\begin{equation}\label{eq:strait_line_distance_objective}
    \Gamma_s(\vect{x}) = \frac{\|\vect{x}-\vect{x}_{h_f}\|_2}{d_{\textit{max}}},
\end{equation}
where $d_{\textit{max}}$ is the maximum distance in the region $\mathcal{D}$ from the goal $\vect{x}_{h_f}$. This portion of the objective function is illustrated in the right panel of Figure~\ref{fig:lpp_objectives}. The heatmap shows the normalized distance to the high-priority agent's goal, guiding the low-priority agents to take more samples in that direction.

The overall objective function is the sum of the three parts:
\begin{equation} \label{eq:objective_function_total}
\Gamma(\vect{x},\mylist{x}_e,\mylist{R})= -\alpha_e\Gamma_e(\vect{x},\mylist{x}_e) + \alpha_u\Gamma_u(\vect{x},\mylist{R}) + \alpha_s\Gamma_s(\vect{x}),    
\end{equation}
where $\alpha_e,\alpha_u,$ and $\alpha_s$ are the weights for the exploration, uncertainty reduction, and goal objective functions, respectively. Because we minimize the objective function, we subtract the exploration objective function, Equation~\eqref{eq:exploration_objective}. This results in maximizing the likelihood of discovering an undiscovered radar.

\subsection{Multi-Agent Optimization}
Given multiple low-priority agents, we require a coordination strategy that prevents agents from converging on identical high-reward locations. To address this, each agent independently selects an optimal waypoint $\vect{x}_{l,i}^w$ (for the $i^{\textit{th}}$ low-priority agent) and travels toward this waypoint via a minimum-time path. However, independent waypoint selection alone is only effective if the chosen waypoints are explicitly deconflicted.

Waypoint deconfliction is achieved by prioritizing waypoint selection based on agent proximity to the radar with the greatest uncertainty, measured by the determinant of its covariance. Waypoint optimization is triggered either when an agent reaches its current waypoint or when a fixed time horizon $T_h$ elapses. Agents sequentially select their waypoints in priority order, minimizing Equation~\eqref{eq:objective_function_total}, and subsequently share their selected waypoints and projected paths to inform the decisions of other agents.

\begin{algorithm}
    \caption{Multi-agent Path Planning for Low Priority Agents}\label{alg:lp_path_planning}
    \begin{algorithmic}[1]
    \STATE{\textbf{input: } $\mylist{R}, \mylist{X}_{e}$}
    \STATE{\textbf{output: } $\vect{x}_{l,i}^w \;\forall i \in \{1,\hdots,N_l\}$}
    \STATE{$t_{h,i} = T_h \; \forall i \in \{1,\hdots,N_l\}$} 
    \WHILE{High priority path not found}
    \IF{$t_{h,i} > T_h \textbf{ for any } i \in \{1,\hdots,N_l\}$}
    \STATE{$j = \algfunc{argmax}_{(\mu_{\vect{x}_{r,j}},\Sigma_{\vect{x}_{r,j}})\in \mylist{R}}\operatorname{det}(\Sigma_{\vect{x}_{r,j}})$}\label{alg:find_most_uncertain}
    \STATE{$\mylist{d}_{i,j} = \{||\vect{x}_{l,i}(t) - \mu_{\vect{x}_{r,j}}||_2\}_{\forall i \in \{1,\hdots,N_l\}}$}\label{alg:find_dist_to_most_uncertain}
    \STATE{$\mylist{i} = \algfunc{argsort}(\mylist{D}_{i,j})$}\label{alg:sort_by_dist_to_most_uncertain}
    \STATE{$\mylist{R}^+\leftarrow \mylist{R},\;\mylist{X}_e^+\leftarrow \mylist{X}_e$}\label{alg:create_temp_lists} 
    \FOR{$i \in \mylist{I}$}
    \STATE{$\vect{x}_{l,i}^w = \algfunc{argmin}_{\vect{x}\in \mathcal{D}}(\Gamma(\vect{x},\mylist{X}_e^+,\mylist{R}^+))$}\label{alg:optimize_waypoint}  
    \STATE{$\mylist{x}_{e,i}^+ = \Bigl\{\vect{x} |\vect{x}=q\frac{\vect{x}_{l,i}^w-\vect{x}_{l,i}(t)}{||\vect{x}_{l,i}^w-\vect{x}_{l,i}(t)||_2}\Delta_{t_e}v_l \;  \forall q \in $} \label{alg:create_future_path} \\
    { \hspace{0.5em} $  \Bigl\{1\hdots,\algfunc{int}\!\!\left(\frac{||\vect{x}_{l,i}^w-\vect{x}_{l,i}(t)||_2}{\Delta_{t_e}v_l}\right)\Bigr\}\Bigr\}$}\COMMENT{Sample future path}
    \STATE{$\mylist{x}_e^+ = \mylist{x}_e^+\bigcup\mylist{x}_{e,i}^+$}\label{alg:combine_future_past_path}
    \STATE{$\mylist{R}^+ = \{(\mu_{\vect{x}_{r,k}},\operatorname{cov}^+(\vect{x}_{l,i}^w,\mu_{\vect{x}_{r,k}},\Sigma_{\vect{x}_{r,k}})\forall k \in \{1,\hdots,N_{\hat{r}}\}\}$}\label{alg:update_future_covariance}
    \STATE{$\algfunc{transmit}(\mylist{R}^+,\mylist{x}_e^+)$}\label{alg:transmit_path_info}
    \STATE{$t_{h,i} = 0$}
    \ENDFOR
    \ENDIF
    \ENDWHILE
    \end{algorithmic}
\end{algorithm}

Algorithm~\ref{alg:lp_path_planning} shows our low-priority path planning scheme. The input to the algorithm is the current estimate of the radar parameters and locations ($\mylist{R}$) and the combined path history of all agents ($\mylist{x}_e$). The output is a waypoint $\vect{x}_{l,i}^w$ for each low-priority agent. If the time since an agent last planned its path ($t_{h,i}$) exceeds the threshold $T_h$, or when any agent reaches their waypoint, new waypoints are found for all agents. The optimization process proceeds as follows:
\begin{itemize}
    \item The radar with the highest uncertainty is identified by maximizing the determinant of its covariance matrix $\Sigma_{\vect{x}_{r,j}}$ (line~\ref{alg:find_most_uncertain}).
    \item Each agent's distance to this radar is computed: $||\vect{x}_{l,i}(t) - \mu_{\vect{x}_{r,j}}||_2$ (line~\ref{alg:find_dist_to_most_uncertain}).
    \item Agents are sorted in increasing order of distance to the most uncertain radar (line~\ref{alg:sort_by_dist_to_most_uncertain}).
    \item Temporary copies of radar estimates $\mylist{R}^+$ and path history $\mylist{x}_e^+$ are initialized (line~\ref{alg:create_temp_lists}).
    \item The $i^\textit{th}$ agent in sorted order optimizes its waypoint $\vect{x}_{l,i}^w$ by minimizing the objective function $\Gamma$ (line~\ref{alg:optimize_waypoint}).
    \item The agent simulates a future path from its current location $\vect{x}_{l,i}(t)$ to the selected waypoint $\vect{x}_{l,i}^w$, using a step size of $\Delta_{t_e}$ (line~\ref{alg:create_future_path}).
    \item This path $\mylist{x}_{e,i}^+$ is appended to the combined future path history $\mylist{x}_e^+$ (line~\ref{alg:combine_future_past_path}).
    \item The radar covariances are updated assuming a measurement at $\vect{x}_{l,i}^w$ (line~\ref{alg:update_future_covariance}).
    \item The updated radar estimates $\mylist{R}^+$ and path history $\mylist{x}_e^+$ are transmitted to other agents (line~\ref{alg:transmit_path_info}).
\end{itemize}

To reduce communication bandwidth, agents may transmit only new or updated information, as described in~\cite{norton2023decentralized}. This process is repeated until all agents have selected deconflicted waypoints.

\section{High-Priority Path Planning} \label{sec:HP_path_planning}
The goal of the high-priority path planner is to find a safe path that starts at $\vect{x}_{h_0}$ and ends at $\vect{x}_{h_f}$, while keeping the PD along the trajectory less than a threshold $P_{D,t}$. The high-priority agent does not have perfect knowledge of the region; it only has estimates of the radar parameters $\mylist{R}$ that the low-priority agents have found. The high-priority agent must account for this uncertainty when planning trajectories. We do this by approximating the probability that the true PD (Equation~\eqref{eq:overall_PD}) is less than a certain threshold (Equation~\eqref{eq:prob_pd_less_than_threshold}). 

In this section, we present two path-planning methods for the high-priority agent, one for the deterministic case, where the radar parameters are all known, and the other for the uncertain case, where the radar parameter estimates and prior beliefs are used. We start by assuming the parameters are known to provide a baseline for comparison and to illustrate the core trajectory optimization approach before extending it to handle uncertainty. Both rely on the use of Voronoi diagrams to find initial feasible trajectories, then using interior point optimization algorithms (IPOPT~\cite{IPOPT}) to optimize a B-spline trajectory which accounts for the path safety (PD) and kinematic feasibility constraints (velocity, curvature, turn rate).

\subsection{Deterministic High-Priority Path Planner}\label{sec:deterministic_HP_planner}
\begin{figure*}[ht]
    \centering
\includegraphics[width=0.95\linewidth,trim={.0cm .1cm .3cm 0.8cm},clip]{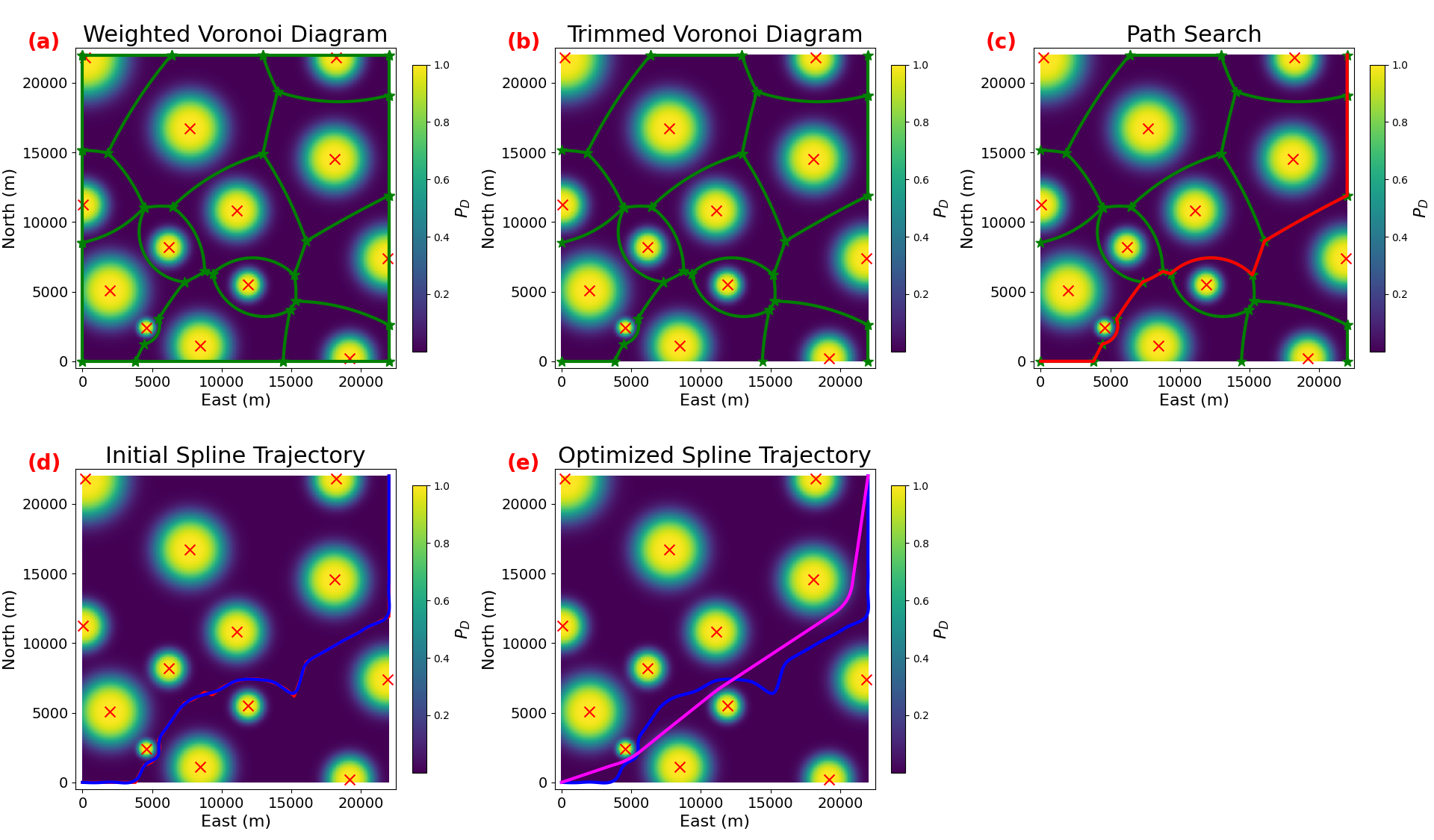}
    \caption{This figure shows the steps of our deterministic high-priority path planning algorithm (Algorithm~\ref{alg:det_high_priority_path_planning}). In the top left figure, the weighted Voronoi diagram is found. Then all edges where $P_D\leq P_{D,t}$ for any point along the edge are removed. In the top right figure, a path search,  A*, is used to find the shortest path through the graph. In the bottom left, a B-spline is fit to the A* path. Then, in the bottom center, the optimized B-spline is shown. All figures show the locations of the radar as red x's with the PD at every location depicted by the heatmap.}\label{fig:det_hpp}
\end{figure*}
For the deterministic case, we assume that the high-priority agent has knowledge of all the radar's parameters, $\vect{x}_{r,j}, P_{T,j}, G_{T,j}, G_{R,j}, \lambda, \tau_{p,j}, T_{s,j}$ and $L_j$. Using this information, we wish to find a feasible trajectory where the PD along the entire trajectory is below a threshold $P_{D,t}$, e.g., $P_{D}(\vect{p}(t)) \leq P_{D,t},$. To do this, we first use a multiplicatively weighted Vornoi diagram along with the A* graph search algorithm~\cite{hart1968formal}. We then fit a B-spline trajectory to the path and use a heuristic to ensure that velocity constraints are met. Finally, we use that B-spline trajectory as a seed trajectory to an interior point optimization algorithm, which finds the minimum time trajectory with PD, velocity, turn rate, and curvature constraints. 

To find an initial feasible trajectory, we use a multiplicatively weighted Voronoi diagram. We show that the minimum PD boundary between two radars is equivalent to finding the edges of a multiplicatively weighted Voronoi cell. We then use the cell edges to find the boundary of equal PD between two radar stations.  For the $j^\textit{th}$ and $k^\textit{th}$ radars, this yields
\begin{equation}
    \operatorname{exp}\left(\frac{\operatorname{ln}P_{fa,j}}{\frac{P_{E,j}G_{R,j}\lambda^2\sigma_i \tau_{p,j}}{(4\pi)^3R_{i,j}^4\kappa T_{s,j}L_j}+1}\right) =\\ \operatorname{exp}\left(\frac{\operatorname{ln}P_{fa,k}}{\frac{P_{E,k}G_{R,k}\lambda^2\sigma_i \tau_{p,k}}{(4\pi)^3R_{i,k}^4\kappa T_{s,k}L_k}+1}\right).
\end{equation}
If we assume that $P_{fa,j} = P_{fa,k}$, finding the boundary where the PD is equal for both radars is the same as finding the region where the SNR for both radars is equal:
\begin{equation}
    {\frac{P_{E,j}G_{R,j}\lambda^2\sigma_i \tau_{p,j}}{(4\pi)^3R_{i,j}^4\kappa T_{s,j}L_j}} = {\frac{P_{E,k}G_{R,k}\lambda^2\sigma_i \tau_{p,k}}{(4\pi)^3R_{i,k}^4\kappa T_{s,k}L_k}}.
\end{equation}
Rearranging this equation to the standard form for multiplicatively weighted Voronoi diagrams (Equation~\eqref{eq:mult_voronoi}), we get
\begin{equation}
    \frac{R_{i,j}}{\left(\frac{\kappa T_{s,j}L_j}{P_{E,j}G_{R,j}\lambda^2\sigma_i \tau_{p,j}}\right)^{\frac{1}{4}}} = \frac{R_{i,k}}{\left(\frac{\kappa T_{s,k}L_k}{P_{E,k}G_{R,k}\lambda^2\sigma_i \tau_{p,k}}\right)^{\frac{1}{4}}}.
\end{equation}
This shows that the Apollonius circle, as defined in weighted Voronoi theory, corresponds to the set of points where the PD is equal between two radar systems.
Consequently, the boundary (i.e., the path of equal PD between two radars) is equivalent to a multiplicatively weighted Voronoi edge, constructed using the radar locations as generator points and the following weight for the $j^\textit{th}$ radar is
\begin{equation}\label{eq:det_weights}
    w_j = \left(\frac{\kappa T_{s,j}L_j}{P_{E,j}G_{R,j}\lambda^2\sigma_i \tau_{p,j}}\right)^{\frac{1}{4}}.
\end{equation}
We denote the set of Voronoi edges generated using the radar locations $\mylist{x}_r = \{\vect{x}_{r,1},\hdots,\vect{x}_{r,N_r}\}$ and the weights $\mylist{w}_r =\{w_1,\hdots,w_{N_r}\}$ found using Equation~\eqref{eq:det_weights} as $\mylist{E}(\mylist{x}_r,\mylist{w}_r)$.

\begin{algorithm}
\caption{Deterministic High-Priority Path Planning}\label{alg:det_high_priority_path_planning}
    \begin{algorithmic}[1]
        \STATE{\textbf{Inputs: $\mylist{x}_r,\mylist{w}_r,\vect{x}_{h_0},\vect{x}_{h_f}$}}
        \STATE{\textbf{Output: $\mylist{C},\vect{t}_k$}}
        \STATE{$\mylist{E} = \algfunc{computeWeightedVoronoiEdges}(\mylist{x}_r,\mylist{W}_r)$}\label{alg:det_compute_weighted_edges}
        \STATE{$\mylist{V} = \algfunc{computeWeightedVoronoiVertices}(\mylist{x}_r,\mylist{W}_r)$}\label{alg:det_compute_weighted_verticies}
        \STATE{$\mylist{V},\mylist{E} = \algfunc{intersectVoronoiEdgesWithBoundary}(\mylist{V},\mylist{E})$}\label{alg:intersect_det_edges_with_boundary}
        \STATE{$\mylist{E} = \algfunc{trimInfeasibleEdges}(\mylist{E})$}\label{alg:det_trim}
        \STATE{$\mylist{V}_\textit{opt},\mylist{E}_\textit{opt} = \algfunc{shortestPathThroughGraph}(\mylist{V},\mylist{E},\vect{x}_{h_0},\vect{x}_{h_f})$}
        \STATE{$\mylist{C}_0 = \algfunc{fitSplineToPath}(\mylist{V}_\textit{opt},\mylist{E}_\textit{opt},N_c,p) $}\label{alg:fit_spline_to_path}
        \STATE{$t_f = \algfunc{checkVelocityConstraint}(\mylist{C}_0,v_\textit{lb},v_\textit{ub})$}
        \STATE{$\mylist{C}_\textit{opt},t_f = \algfunc{optimizeTrajectory}(\mylist{C}_0,t_f,\vect{x}_{h_0},\vect{x}_{h_f})$}
    \end{algorithmic}
\end{algorithm}

Algorithm~\ref{alg:det_high_priority_path_planning} presents our high-priority path planning method. It takes as input the radar locations and their associated weights, which are computed using Equation~\eqref{eq:det_weights}. The algorithm outputs a set of control points, $\mylist{C}_\textit{opt}$, along with knot points (determined using the final time $t_f$) that define an optimized B-spline trajectory. The steps of this algorithm are shown in Figure~\ref{fig:det_hpp}.

\emph{ Computing Weighted Voronoi Diagrams.}
The algorithm's first step is to compute the weighted Voronoi diagram, specifically its edges $\mylist{E}$ and vertices $\mylist{V}$, using the method outlined in~\cite{held2020efficient} with weights given by Equation~\eqref{eq:det_weights}. These diagrams represent proximity relationships among the radar nodes, factoring in their weights. The edges in this case are arcs of circles defined from the weighted Vornoi diagram. This step is shown in lines~\ref{alg:det_compute_weighted_edges} and~\ref{alg:det_compute_weighted_verticies}. Figure~\ref{fig:det_hpp} (a) shows the edges and vertices in green.

\emph{ Intersecting Voronoi Edges with the Boundary.}
In line~\ref{alg:intersect_det_edges_with_boundary}, we intersect the Voronoi edges with the boundary of the operating region $\mathcal{D}$. For a rectangular region defined by lower and upper bounds $\vect{x}_\textit{lb}$ and $\vect{x}_\textit{ub}$, we check each edge in $\mylist{E}$ for intersections with the four sides of the region. Intersections result in new vertices added to $\mylist{V}$, and the corresponding edges are trimmed to lie within the region. Additionally, edges between boundary-intersecting vertices are added to $\mylist{E}$ to preserve connectivity.
This creates a graph where the connections are either arcs of circles defined from the weighted Voronoi diagram or straight lines from the boundary. The intersected graph is shown in Figure~\ref{fig:det_hpp} (a). 

\emph{ Trimming Edges Based on PD Constraints.}
In line~\ref{alg:det_trim}, we remove the edges in $\mylist{E}$ whose maximum PD along the edge exceeds a threshold $P_{D,t}$. The maximum PD along an edge occurs at the point where the weighted distance to the generator points is minimized. Each edge is an arc of an Apollonius circle, defined by a center $\vect{x}_{e,i}$ (Equation~\eqref{eq:apollonius_center}), radius $r_{e,i}$ (Equation~\eqref{eq:apollonius_radius}), and start/stop angles $\theta_{0,e,i}$ and $\theta_{f,e,i}$. For each edge, we identify the point on the arc that is closest in Euclidean distance to the generator points. Because the arc satisfies a constant ratio of distances to the two generators, the point that minimizes the Euclidean distance to one generator also minimizes the weighted distance to both generator points, and thus corresponds to the location of maximum PD. This point is used to determine whether the edge should be trimmed from the diagram. The closest point on an arc to a point is found using the following equations:
\begin{align}
\theta_{i,j} &= \arctan \left( \frac{y_{r,j} - y_{e,i}}{x_{r,j} - x_{e,i}} \right)\notag\\
\vect{x}_{i,j} \!\!&=\!\! \begin{cases}
    \vect{x}_{e,i} + r_{e,i}\begin{bmatrix}
        \cos{\theta_{i,j}}\\
        \sin{\theta_{i,j}}
    \end{bmatrix}\!\!,  \theta_{0,e,i}\leq \theta_{i,j} \leq \theta_{f,e,i}\\
    \vect{x}_{e,i} + r_{e,i}\begin{bmatrix}
        \cos{\theta_{0,e,i}}\\
        \sin{\theta_{0,e,i}}
    \end{bmatrix}\!\!,  |\theta_{0,e,i} - \theta_{i,j}|\leq |\theta_{f,e,i} - \theta_{i,j}|\\
    \vect{x}_{e,i} + r_{e,i}\begin{bmatrix}
        \cos{\theta_{f,e,i}}\\
        \sin{\theta_{f,e,i}}
    \end{bmatrix}\!\!,  |\theta_{f,e,i} - \theta_{i,j}|\leq |\theta_{f,e,i} - \theta_{i,j}|
\end{cases}
\end{align}

The PD at this closest point is found using Equation~\eqref{eq:overall_PD}. If the PD is greater than the threshold $P_{D,t}$, that edge is removed from $\mylist{E}.$
The trimmed graph is shown in Figure~\ref{fig:det_hpp} (b).

\emph{ Graph Construction and Path Search.}
Using the trimmed set of vertices $\mylist{V}$ and edges $\mylist{E}$ a graph is created. We then use the A* algorithm~\cite{hart1968formal} to find the shortest path through the graph, where the distance between the nodes is the arc distance of the edge connecting the nodes, starting at the node corresponding to $\vect{x}_{h_0}$ and ending at $\vect{x}_{h_f}$. This path is defined as an ordered list of vertices $\mylist{V}_\textit{opt}$ and edges $\mylist{E}_\textit{opt}$. This shortest path is shown in red in Figure~\ref{fig:det_hpp} (c).

\emph{ Spline Fitting and Time Adjustment.}
The next step is to fit a B-spline to the shortest path (line~\ref{alg:fit_spline_to_path}). The path is first sampled at a spatial frequency of $\Delta_{x}$, producing a discretized path represented as a sequence of points. We apply the least-squares algorithm proposed in~\cite{dierckx1995curve} and implemented in the Scipy library~\cite{2020SciPy-NMeth} to fit a B-spline with $N_c$ control points, degree $p$, and internal knots spaced over the interval $[0,1]$ to the discrete path. To satisfy velocity constraints, we use a heuristic approach: we evaluate the spline’s velocity using Equation~\eqref{eq:velocity_constraint} and identify its maximum value. If the constraint is violated, $t_f$ is increased, knot points are recomputed over $[0, t_f]$, and the process is repeated until the constraint is met. The spline fit to the shortest path is shown in blue in Figure~\ref{fig:det_hpp}(d). This path becomes the initial guess to a B-spline optimization problem.

\emph{ Trajectory Optimization.}
Finally, an interior-point optimization algorithm (IPOPT~\cite{IPOPT}) refines the trajectory to minimize travel time, while conforming to physical feasibility constraints. The optimization problem is:

\begin{subequations} \label{eq:optimization}
\begin{align}
\mylist{C}_{opt},t_f =\operatorname*{argmin}_{\mylist{C}, t_f} \quad& t_f \label{eq:optimization_util_det} \\
\text{s.t.} \quad&
\vect{p}(0) = \vect{x}_{h_0} \label{eq:start_constraint_det} \\
&\vect{p}(t_f) = \vect{x}_{h_f} \label{eq:end_constraint_det} \\
&\vect{p}(\vect{t}_s) \in \mathcal{D} \label{eq:optimization_path_boundary_constraint_det} \\
&v_{lb} \leq v(\vect{t}_s) \leq v_{ub} \label{eq:velocity_constraint_det} \\
&u_{lb} \leq u_A(\vect{t}_s) \leq u_{ub} \label{eq:turn_rate_constraint_det} \\
&-\kappa_{ub} \leq \kappa_A(\vect{t}_s) \leq \kappa_{ub} \label{eq:curve_constraint_det} \\
&P_{D}(\vect{p}(\vect{t}_s)) \leq P_{D,t} \label{eq:det_pd_constraint},
\end{align}
\end{subequations}
where $\vect{t}_s = \{0, \Delta_s, 2\Delta_s,\hdots, t_f\}$ are the discrete points in time where the constraints are evaluated, $\Delta_s = t_f/N_s$ is the time spacing between constraint samples and $N_s$ is the number of constraint samples. We discretely sample the constraints, which means that the constraints could be violated between samples. However, the likelihood of this is reduced with an increasing number of samples.  

The objective of Equation~\eqref{eq:optimization_util_det}, is to find the minimum time trajectory, starting at $\vect{x}_{h_0}$ and ending at $\vect{x}_{h_f}.$ The constraint in Equation~\eqref{eq:optimization_path_boundary_constraint_det} ensures that the agent stays within the operating region $\mathcal{D}$. The next three constraints, Equations~\eqref{eq:velocity_constraint_det} through Equation~\eqref{eq:curve_constraint_det}, show the constraints of kinematic feasibility, velocity, turn rate, and curvature. These ensure that the agent can physically follow the planned trajectory. The final constraint, Equation~\eqref{eq:det_pd_constraint}, ensures that the PD of the agent (found using Equation~\eqref{eq:overall_PD}) is below the threshold. The magenta line in Figure~\ref{fig:det_hpp} (e) shows the minimum time spline trajectory found.

\subsection{Uncertain High-Priority Path Planner}\label{sec:uncertain_HP_planner}
\begin{figure*}[ht]
    \centering
\includegraphics[width=0.95\linewidth,trim={.0cm 0.1cm 0.1cm 0.8cm},clip]{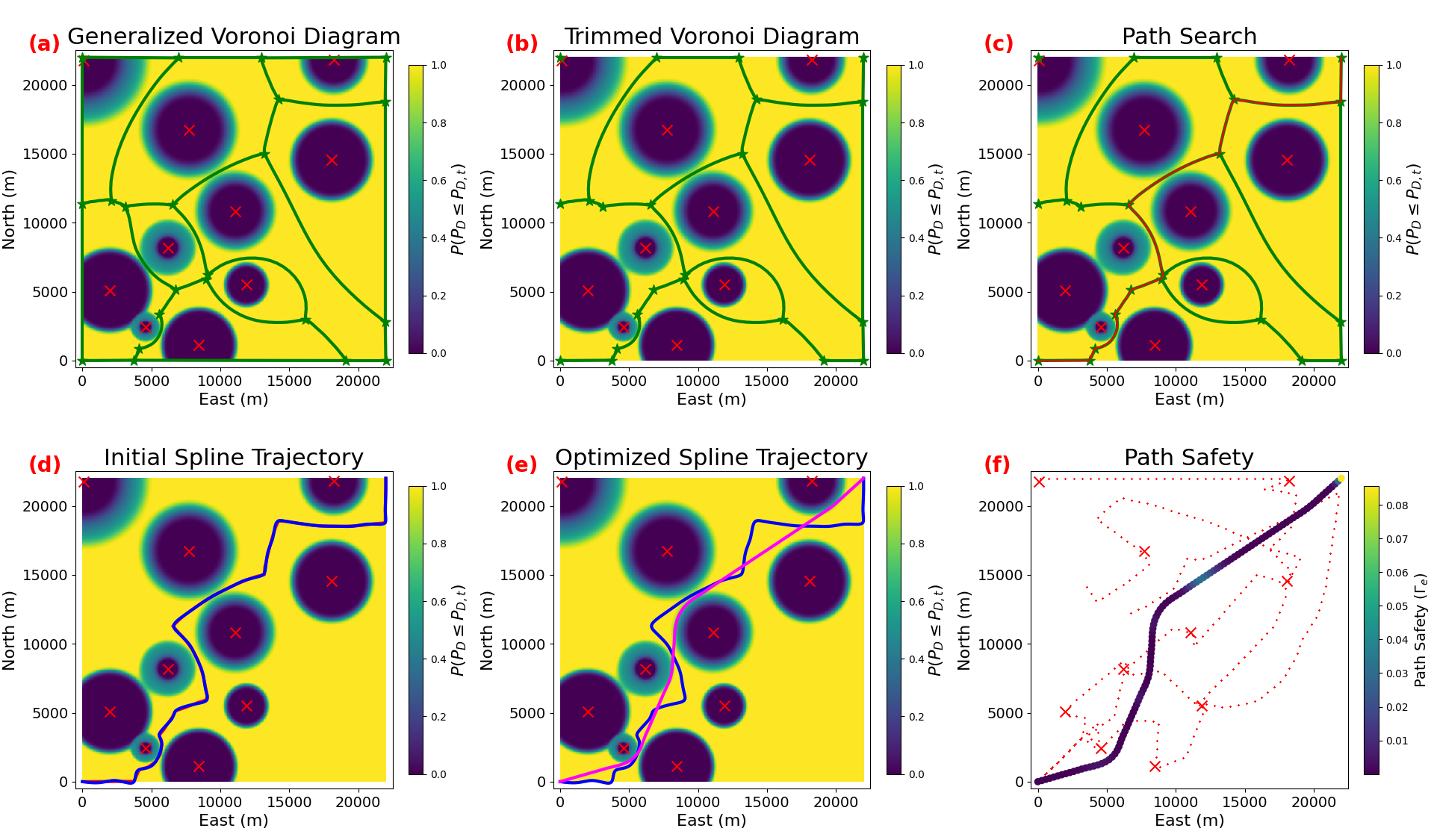}
    \caption{This figure shows the steps of our uncertain high-priority path planning algorithm (Algorithm~\ref{alg:uncertain_high_priority_path_planning}). In the top left figure, the generalized Voronoi diagram is found using $P(P_D\leq P_{D,t})$ as the criterion. Then all edges where $P(P_D\leq P_{D,t})<\epsilon$ are removed. In the top right figure, A* is used to find the shortest path through the graph. In the bottom left, a B-spline is fit to the A* path. Then, in the bottom center, the optimized B-spline is shown. Finally, the bottom left image shows the path safety metric. In this figure, the low-priority agents' paths are shown in red. All figures show the current estimated location of the radar as red x's. The first five show $P(P_D\leq P_{D,t})$ as the heatmap.}\label{fig:uncertain_path_planning}
\end{figure*}

In the deterministic setting, the high-priority path planner assumes full knowledge of all radar parameters. Under this assumption, a multiplicatively weighted Voronoi diagram is used to define regions of equal PD between radar systems. This Voronoi-based representation helps generate an initial feasible trajectory, which is then refined using an IPOPT-based optimization that directly enforces a constraint on PD: the agent’s trajectory must remain below a fixed detection threshold, $P_{D,t}$, at all points along the path.

When radar parameters are uncertain, however, the high-priority agent must account for both estimated values and prior distributions over unknown parameters. In this setting, we extend the deterministic approach by introducing a probabilistic safety constraint: rather than requiring the PD to stay strictly below a threshold, we require the probability that the PD remains below that threshold to exceed a specified confidence level, $\left(P(P_D\leq P_{D,t})\geq \epsilon\right)$. This formulation ensures that the planned path remains sufficiently safe, even in the presence of parameter uncertainty.

The use of Voronoi diagrams also changes under uncertainty. Instead of using a multiplicatively weighted Voronoi diagram, we construct a generalized Voronoi diagram based on the probability that the PD at a location is below the threshold. Each point in the operational domain is assigned to the radar system with the highest likelihood of violating the safety constraint, and boundaries are defined accordingly. This risk-aware partitioning guides the initial trajectory generation, which is then refined using the same B-spline-based optimization framework, but with probabilistic constraints replacing deterministic ones.

The remainder of the algorithm proceeds in the same way as in the deterministic case, with the exception of the path safety constraint used.
\begin{algorithm}
\caption{Uncertain High-Priority Path Planning}\label{alg:uncertain_high_priority_path_planning}
    \begin{algorithmic}[1]
        \STATE{\textbf{Inputs: $\mylist{R}_t, \mylist{x}_{e,t}, \vect{x}_{h_0}, \vect{x}_{h_f}$}}
        \STATE{\textbf{Output: $\mylist{C}, \vect{t}_k$}}
        \STATE{$\mylist{E}, \!\mylist{V}\!\! =\!\! \algfunc{computeGeneralizedVoronoiEdgesAndVertices}(\mylist{R})$}
        \STATE{$\mylist{E} = \algfunc{trimInfeasibleEdges}(\mylist{E},\mylist{R})$}\label{alg:uncertain_trim}
        \STATE{$\mylist{V}_\textit{opt},\mylist{E}_\textit{opt} = \algfunc{shortestPathThroughGraph}(\mylist{V},\mylist{E},\vect{x}_{h_0},\vect{x}_{h_f})$}
        \STATE{$\mylist{C}_0 = \algfunc{fitSplineToPath}(\mylist{V}_\textit{opt},\mylist{E}_\textit{opt},N_c,p) $}\label{alg:_uncertain_fit_spline_to_path}
        \STATE{$t_f = \algfunc{checkVelocityConstraint}(\mylist{C}_0,v_\textit{lb},v_\textit{ub})$}
        \STATE{$\mylist{C}_\textit{opt},t_f = \algfunc{optimizeTrajectory}(\mylist{C}_0,t_f,\vect{x}_{h_0},\vect{x}_{h_f},\mylist{R})$}
        \STATE{$P_{\textit{max}} = \algfunc{max}_{t_s\in\vect{t}_s}\Gamma_e(\vect{p}(t_s),\mylist{x}_{e,t})$}
        \IF{$P_{\textit{max}}\leq P_s$}
        \STATE{Dispatch high-priority agent}
        \ENDIF
    \end{algorithmic}
\end{algorithm}

Algorithm~\ref{alg:uncertain_high_priority_path_planning} outlines our method. The input of the algorithm is the current estimate of the radar parameters and the covariances $\mylist{r}_t,$ given from the EKF described in Section~\ref{sec:radar_localization}, the current combined agent path history list $\mylist{x}_{e,t}$, at time $t$,  and the start $\vect{x}_{h_0}$ and destination $\vect{x}_{h_f}$ locations of the high priority agent. 

The first step in the algorithm is to find the generalized Voronoi edges and vertices using $P(P_D \leq P_{D,t})$ as the criterion. This is accomplished using the grid-based algorithm shown in Algorithm~\ref{alg:grid_voronoi}.
\begin{algorithm}
\caption{Grid Based Generalized Voronoi Diagram}\label{alg:grid_voronoi}
    \begin{algorithmic}[1]
        \STATE{\textbf{Inputs: $\mylist{R}$}}
        \STATE{\textbf{Output: $\mylist{E}, \mylist{V}$}}
        \STATE{$\mylist{X}_t = \algfunc{evenlySpaceSampeles}(\mathcal{D})$}
        \STATE{$\mylist{A} = \{\algfunc{argmin}_j(P(P_D(\vect{x}_{t,i},\theta_{u,j},\theta_{e,j})\leq P_{D,t})\}\forall \vect{x}_{t,i}\in \mylist{x}_t$}\label{alg:grid_vor_cell_assignment}
        \STATE{$\mylist{N} = \algfunc{findAdjecentCellsUsingVoronoi}(\mylist{R})$}
        \STATE{$\mylist{E},\mylist{V} = \algfunc{findGeneralizedVoronoiEdges}(\mylist{N},\mylist{A})$}
    \end{algorithmic}
\end{algorithm}
We begin by generating an evenly spaced set of test points in the region $\mylist{x}_t$. At each test point, we evaluate the likelihood that the true PD falls below the threshold for each radar.  Each point is then assigned to the Voronoi cell corresponding to the radar that minimizes this likelihood, as shown in line~\ref{alg:grid_vor_cell_assignment}. 

This process provides the generalized Voronoi cells; however, we need to find the edges and vertices between these cells. To find the edges, we look at the contours around each cell. For each pair of neighboring cells, we find the portion of each cell's contours that overlap. This overlapped portion of the cell boundary is the generalized Voronoi edge between the cells. We then fit a third-order B-spline to this boundary, resulting in each edge being defined using control points and knot points $e_{i,j} = (\mylist{C}_{ij},\vect{t}_{i,j})$. The Voronoi vertices are the points where the generalized Voronoi edges intersect.

After finding the generalized Voronoi edges and vertices, Algorithm~\ref{alg:uncertain_high_priority_path_planning} proceeds similarly to Algorithm~\ref{alg:det_high_priority_path_planning}. Infeasible edges are trimmed in Line~\ref{alg:uncertain_trim}. Because there is no closed-form solution to the location of the lowest $P(P_D\leq P_{D,t})$, we sample the edge and check many different locations. An edge is removed if the probability that PD is below the threshold is less than $\epsilon$ at any of the sampled points. After infeasible edges are removed, we use the A-Star algorithm to find the shortest path through the graph created using the generalized Voronoi vertices and trimmed edges. We then sample the path and fit a B-spline to the path using $N_c$ control points and internal knot points defined on $[0,1]$. We use the same heuristic method described for the deterministic path planner to ensure the velocity constraint is met. After a feasible initial trajectory is found, we use IPOPT to refine the trajectory according to the following optimization problem:
\begin{subequations} \label{eq:hpp_opt_uncertain}
\begin{align}
\mylist{C}_{opt},t_f = \operatorname*{argmin}_{\mylist{C},t_f}\quad&t_f\label{eq:uncertain_optimization_util}\\
\text{s.t.}\quad&
 \vect{p}(0) = \vect{x}_{h_0} \label{eq:start_constraint_uncertain}\\
 &\vect{p}(t_f) = \vect{x}_{h_f}\label{eq:end_constraint_uncertain}\\
 &\vect{p}(\vect{t}_s) \in \mathcal{D}\label{eq:optimization_path_boundary_constraint_uncertain}\\
 &v_{lb} \leq v(\vect{t}_s) \leq v_{ub}\label{eq:velocity_constraint_uncertain}\\
&u_{lb} \leq u_A(\vect{t}_s) \leq u_{ub}\label{eq:turn_rate_constraint_uncertain}\\
&-\kappa_{ub} \leq \kappa_A(\vect{t}_s) \leq \kappa_{ub} \label{eq:curve_constraint_uncertain}\\
&P(P_D(\vect{p}(\vect{t}_s),\mylist{\theta}_{u},\mylist{\theta}_{e})\leq P_{D,t})\geq \epsilon \label{eq:uncertain_pd_constraint}.
\end{align}
\end{subequations}
This optimization problem is the same as the deterministic case (Equation~\eqref{eq:optimization}) except for the final constraint~\eqref{eq:uncertain_pd_constraint}. Instead of using the ground-truth PD evaluated at points along the trajectory, we use the approximate probability that the ground-truth PD is less than $P_{D,t}$. This ensures the probability of the path being safe, $P(P_D \leq P_{D,t}$), is greater than $\epsilon$, accounting for the uncertainty from the estimates and prior beliefs. 

The final step of the algorithm is to test the path to see if it is likely that there are undiscovered radar stations near the path. So far, the algorithm has only accounted for discovered radar, and it has no sense of the danger from undiscovered radar. Using Equation~\eqref{eq:exploration_objective}, we find the probability that an undiscovered radar is on the optimized trajectory. If this probability is greater than a threshold $P_s$ at any point along the trajectory, we determine that the path is not safe enough for the high-priority agent to traverse. If the probability is lower than the threshold at all points on the trajectory, the path is determined to be safe, and the high-priority agent is dispatched. Each step of the algorithm is outlined in Figure~\ref{fig:uncertain_path_planning}.

\section{Results} \label{sec:results}
In this section, we present simulation results for our low-priority and high-priority path planning algorithms. We first present results for the low-priority path planner, showing how it performs with various parameter values and how it compares to a baseline, ``lawnmower" path planner. We next show results for the high-priority path planner.
To test the algorithm in varying situations, we randomly placed 13 different radars in the region with random parameters. We generated 50 different random configurations of the radar using the parameters shown in Table~\ref{tab:parameter_table}. For each random configuration, we sample the radar parameters from a uniform distribution defined by the range provided in the table. For example, the output power is sampled from a uniform distribution from $P_{_t,l}$ to $P_{_t,u}$. Additional parameters used in the simulation experiments are shown in Table~\ref{tab:parameter_table}. 
\begin{table}[h!]
\resizebox{0.95\columnwidth}{!}{
\begin{tabulary}{0.3pt}{| l | c | c | c |}
 \hline
 \textbf{Parameter} & \textbf{Symbol} & \textbf{Algorithm} & \textbf{Value} \\
 \hline
 \hline
 \text{Operational Region Bounds} &  - & - & $(22000,22000)$\\
 \hline
 \text{Radar Output Power Upper} &  $P_{T,u}$ & - & $20000$ Watts\\
 \hline
 \text{Radar Output Power Lower} &  $P_{T,l}$ & - & $0$ Watts\\
 \hline
 \text{Radar Transmit Gain Upper} &  $G_{T,u}$ & - & $20$dB\\
 \hline
 \text{Radar Transmit Gain Lower} &  $G_{T,l}$ & - & $0$dB\\
 \hline
 \text{Radar Receive Gain} &  $G_{R}$ & - & $10$dB\\
 \hline
 \text{Path Loss} &  $L$ & - & $0$dB\\
 \hline
 \text{Radar Wavelength} &  $\lambda$ & - & $99.9$ millimeters\\
 \hline
 \text{Radar Pulse Width} &  $\tau_{p}$ & - & $0.000011$ seconds\\
 \hline
 \text{Radar System Temperature} &  $T_{s}$ & - & $745$ Kelvin\\
 \hline
 \text{Radar Probability of False Alarm} &  $P_{f_a}$ & - & $745$ Kelvin\\
 \hline
 \text{Agent Intercept Antenna Gain} &  $G_{I}$ & - & $1$dB\\
 \hline
 \text{Agent Radar Cross Section} &  $G_{I}$ & - & $0.1$ meters squared\\
 \hline
 \text{Angle of Arrival Measurement Noise Standard Deviation} &  $\sigma_\phi$ & - & $2$ Degrees\\
 \hline
 \text{Received Power Measurement Noise Standard Deviation} &  $\sigma_{S_E}$ & - & $1.0$ microwatt\\
 \hline
 \text{Probability of Intercept Loss} &  $\delta$ & Low Priority & $1000$ \\
 \hline
 \text{Covariance Objective Function Multiple} &  $d_\textit{cov}$ & Low Priority & $1e14$ \\
 \hline
 \text{Re-plan Horizon} &  $T_h$ & Low Priority & $20$ Seconds \\
 \hline
 \text{Agent Path History Period} &  $\Delta_{t_e}$ & Low Priority & $5$ Seconds \\
 \hline
 \text{PD Threshold} &  $P_{D,t}$ & High Priority & $15$\% \\
 \hline
 \text{B-spline Degree} &  $\rho$ & High Priority & $3$ \\
 \hline
 \text{Number of Control Points} &  $N_c$ & High Priority & $40$ \\
 \hline
 \text{Start Location} &  $\vect{x}_{h_0}$ & High Priority & $(0,0)$ \\
 \hline
 \text{Goal Location} &  $\vect{x}_{h_f}$ & High Priority & $(22000,22000)$ \\
 \hline
 \text{Velocity Bounds} &  $v_\textit{lb},v_\textit{ub}$ & High Priority & $100,134$ m/s\\
 \hline
 \text{Turn Rate Bounds} &  $u_\textit{lb},u_\textit{ub}$ & High Priority & $-5,5$ rad/s\\
 \hline
 \text{Maximum Curvature} &  $\kappa_\textit{ub}$ & High Priority & $0.1$ rad/m\\
 \hline
 \text{Probability of PD Threshold} &  $\epsilon$ & High Priority & $0.9$\\
 \hline
 
 \hline
\end{tabulary}}
\caption{\label{tab:parameter_table}Simulation Parameters}
\end{table}

\begin{figure}[ht]
    \centering
    \includegraphics[width=0.95\linewidth,trim={6.0cm 2.9cm 6.9cm 0.3cm},clip]{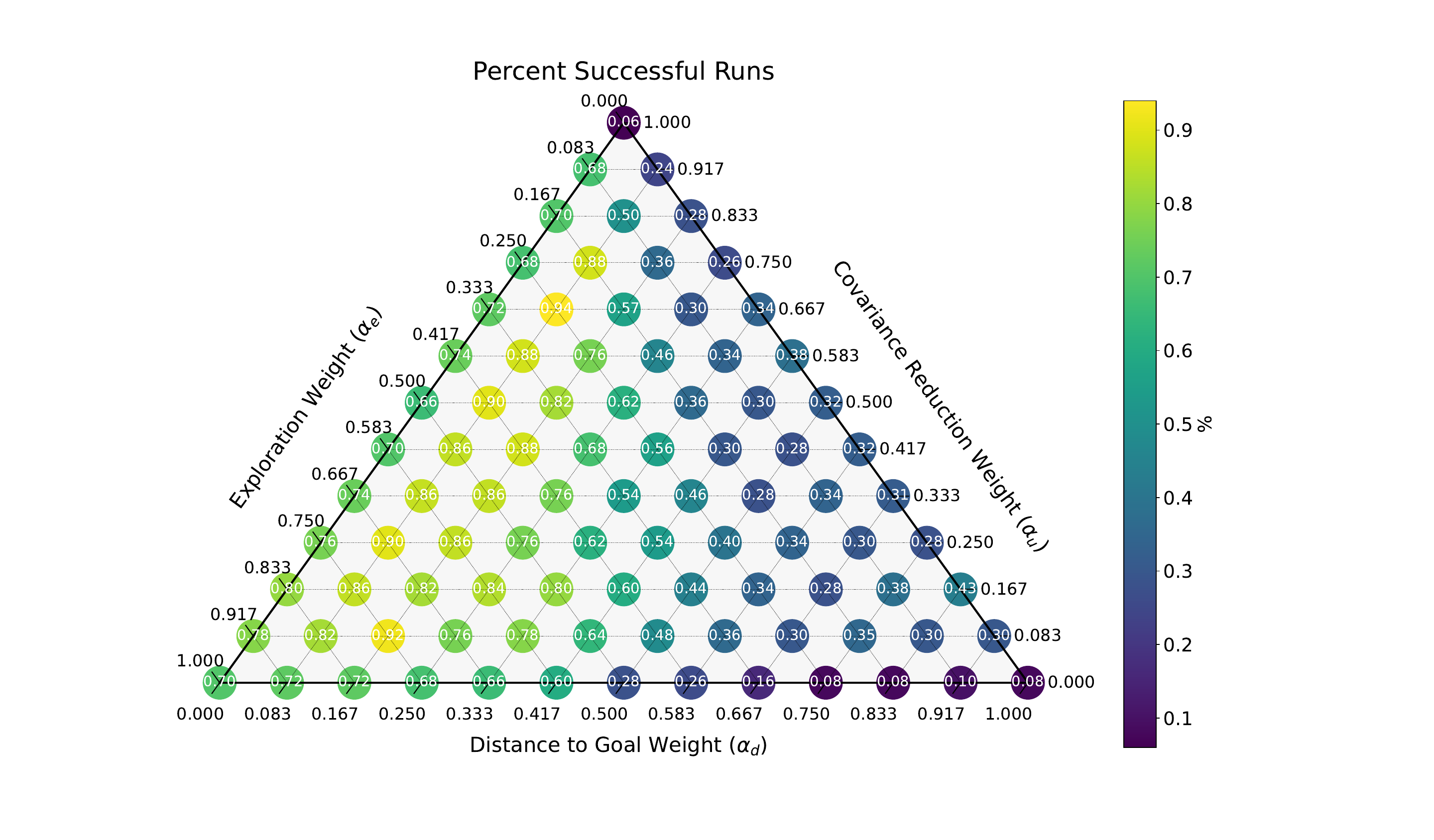}    \caption{This figure shows the percent of random runs where the low-priority agents were able to find the path for the high-priority agents with varying values of the exploration weight $\alpha_e$, the covariance reduction weight $\alpha_u$ and the distance to goal weight $\alpha_s$. The plot is a ternary plot, where each parameter is represented by one side of the triangle. The color of the points corresponds to the percent of runs that were successful.}
    \label{fig:alpha_test_percent}
\end{figure}

\begin{figure}[ht]
    \centering
\includegraphics[width=0.95\linewidth,trim={6.0cm 2.9cm 6.9cm 0.3cm},clip]{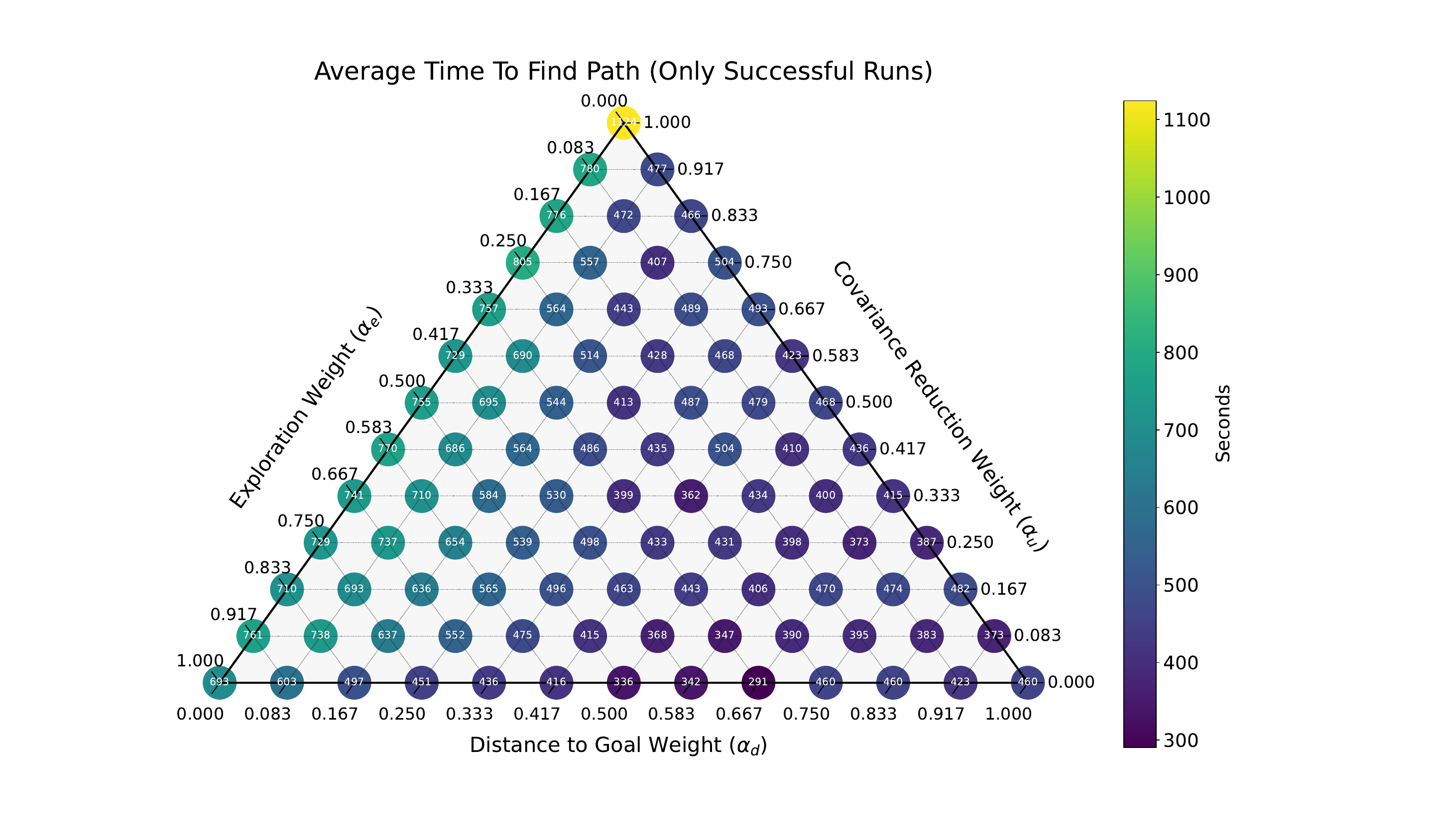}    \caption{This figure shows the average time it took for the low-priority agents to find a path for the high-priority agent with varying levels of parameters. The average time only includes runs that were successful (see Figure~\ref{fig:alpha_test_percent} for how many runs were successful). The color of the points shows the average time it took for the agents to find the path.}
    \label{fig:alpha_test_time}
\end{figure}

We first show how varying weights of the low-priority path planning algorithm's objective function ($\alpha_e,\alpha_u,\alpha_s$) affects performance. We constrain the sum of the tuning parameters to be one, $\alpha_e+\alpha_u+\alpha_s=1$. We then vary the value of each tuning parameter with this constraint and test the algorithm's performance on the 50 random radar configurations. The results of these tests are shown in Figures~\ref{fig:alpha_test_percent} and~\ref{fig:alpha_test_time}. A ternary plot shows the effects of the three parameters. Each side of the triangle represents the value of one parameter being weighted exclusively, with no weight to the other two parameters.
The lines in the plot show the grid where that parameter is the same. 

Figure~\ref{fig:alpha_test_percent} shows the percentage of the random runs where the low-priority agents were able to find a path for the high-priority agent. As seen in the figure, the parameters need to be balanced. At the corners, where a single weight dominates, the success rate is low. Using the figure, we can approximate the ideal weights. The distance from the target weight should be lower, with a value of $\alpha_s=0.083$ showing the best performance. There also needs to be a balance between exploration and covariance reduction, with the best performance occurring at $\alpha_u = 0.667$ and $\alpha_e=0.25$. In this case, the agents were able to find the path for the high-priority agent in $94\%$ of the random runs.  

Figure~\ref{fig:alpha_test_time} shows the average time it took for low-priority agents to find a path for the high-priority agent with varying parameter values. This only includes the successful runs (because in the runs that were not successful, the agents never found a path for the high-priority agent). It can be seen that increasing the distance to the goal weight $\alpha_s$ decreases the average time to find the path; however, in these cases, a smaller percentage of runs were successful. 

\begin{figure}[ht]
    \centering
\includegraphics[width=0.95\linewidth,trim={0.0cm 0.9cm 0.9cm 0.3cm},clip]{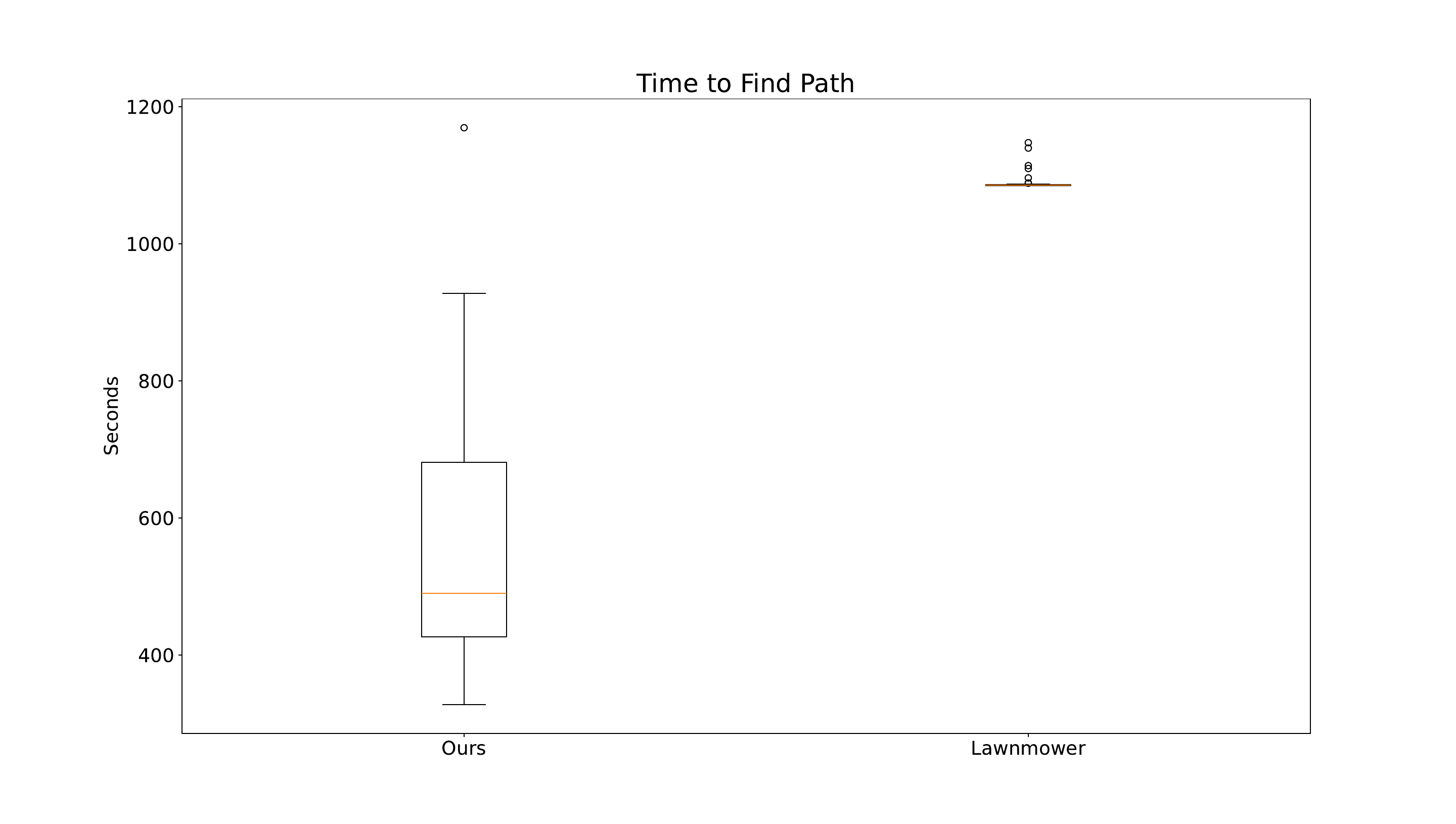}    \caption{This figure shows a box plot distribution of the time it took the low-proirty agents to find the high priority path between the baseline ``lawnmower" path planner and our algorithm using the best parameters: $\alpha_s=0.083$, $\alpha_u = 0.667$, $\alpha_e=0.25.$}
    \label{fig:timing_comparison}
\end{figure}
We next show a comparison of our low-priority path planning algorithm with a baseline, ``lawnmower" path planner. The ``lawnmower" path planner divides the region equally between all agents and creates a sweeping pattern back and forth for each agent in their section. We choose the ``rung" size of the pattern using the exploration objective function Equation~\eqref{eq:exploration_objective}, by picking a threshold of how likely an undiscovered radar would lie in between rungs and solving for a distance. We also allowed the ``lawnmower" path planner to travel back down through the region, splitting the original rungs in half. The agents were limited to running for 1200 seconds for both our algorithm and the ``lawnmower" algorithm. Using our path planner, the low-priority agents were able to find the high-priority path in $94\%$ of the random runs. The ``lawnmower" low-priority agents were able to find the high-priority path in $86\%$ of the runs. The timing results are shown in Figure~\ref{fig:timing_comparison}. As can be seen in the figure, the ``lawnmower" path planner finds the path in roughly the same amount of time each run. This is determined by how long it takes the agents to complete their coverage paths. Our algorithm's times are more distributed because of how its performance depends on the radar layout. Our algorithm is able to find the high-priority path faster in all but one case.

\begin{figure}[ht]
    \centering
\includegraphics[width=0.95\linewidth,trim={0.0cm 0.9cm 0.9cm 0.3cm},clip]{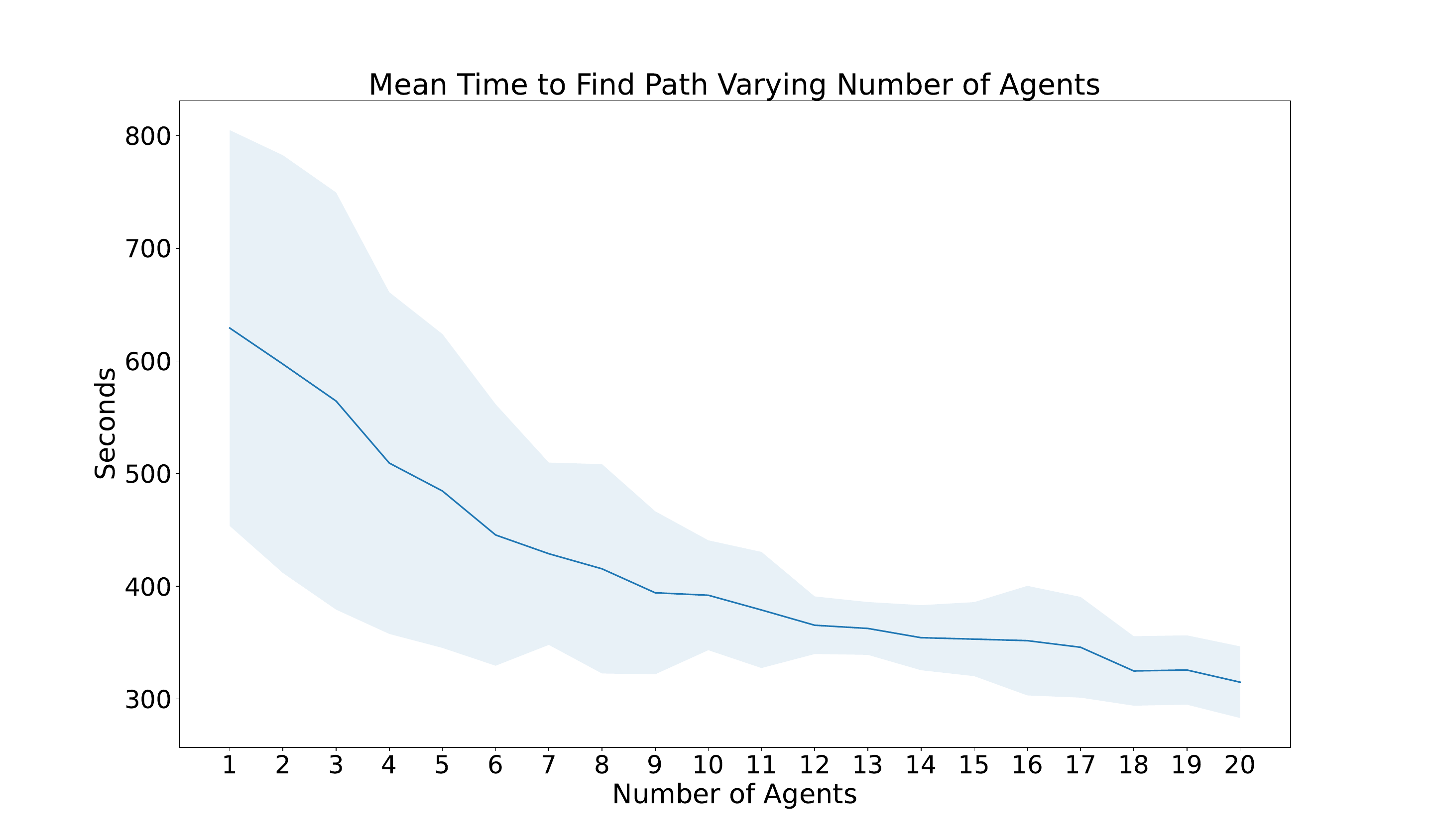}    \caption{This figures shows the average time it took for the low-proirty agents to find a safe path for the high-priority agent using the best parameters from the previous test: $\alpha_s=0.083$, $\alpha_u = 0.667$, $\alpha_e=0.25.$ The time decreases with increasing number of agents because there are more agents to rapidly explore the environment. One standard deviation above and below the mean time are shaded in the figure.}
    \label{fig:numAgents_vs_time}
\end{figure}

Figure~\ref{fig:numAgents_vs_time} shows how our algorithm performs with varying numbers of agents. The plot shows the mean time it took the low-priority agents to find the high-priority path with varying numbers of agents. The same 50 random scenarios were used in this test. The best parameters from the tuning tests were used: $\alpha_s=0.083$, $\alpha_u = 0.667$, $\alpha_e=0.25.$ As expected, with increasing numbers of low-priority agents, the time it takes them to find the high-priority path decreases because there are more agents to explore the environment. 

To show how well our probabilistic high-priority path planning algorithm performs, we provide statistics from the run with 20 low-priority agents and parameters $\alpha_s=0.083$, $\alpha_u = 0.667$, $\alpha_e=0.25$ where the agents were able to find a safe path in all 50 cases. To show how the path planner performs, we first report the average maximum ground truth PD accounted for along the trajectories. We evaluate the ground truth PD along each of the 50 planned paths and find the maximum, then compute the average, which is $11.59\%.$ This is below the $15\%$ threshold used in the algorithm. The maximum ground truth PD was above the $15\%$ threshold in 3 cases, meaning the success rate of the probabilistic threshold was $94\%$. We used a $90\%$ threshold, so we would expect the probabilistic threshold to work in about $90\%$ of the cases if the linearization of the PD was accurate. Our results match our expectations and illustrate the efficacy of our approach in finding probabilistically safe paths.

\section{Conclusion} \label{sec:conclusion}
In this paper, we presented three probabilistic path planning algorithms for operating in contested enemy environments, leveraging a risk-aware framework involving low-cost scouting agents and a high-priority agent tasked with completing a critical mission. The enemy environment was characterized by enemy radar, and we modeled the PD from that radar. We assumed low-priority agents could intercept enemy radar signals and estimate the parameters of enemy radar from the intercepted signals. The low-priority agents were tasked with exploring the environment to find a safe path for the high-priority agent. 

We developed a decentralized path planning algorithm for the low-priority agents that planned paths to locate unknown enemy radar stations and visit informative measurement locations to reduce uncertainty in radar parameter estimates. We presented simulation results comparing our method with a ``lawnmower" baseline and demonstrated improved performance. 

We also presented two algorithms for the high-priority agent: one for the deterministic case, where all radar parameters are known, and the other for the uncertain case, where radar parameters must be inferred. Both approaches rely on a Voronoi diagram and the A* algorithm to generate an initial feasible trajectory through the radar field, which is then refined using an interior point optimization algorithm that accounts for kinematic feasibility and PD constraints. Simulation results demonstrated that the algorithm effectively preserved path safety while incorporating radar parameter uncertainty.

Future work includes accounting for agent RCS that varies with the view angle of the vehicle. In this work, we assumed a constant RCS across all view angles. This extension could be incorporated by using the maximum RCS when selecting the initial trajectory, and allowing the optimization algorithm to adjust the vehicle's orientation to exploit the true RCS variation for improved performance.

%







\bibliography{bib}

\begin{thebibliography}{41}
\newcommand{\enquote}[1]{``#1''}
\providecommand{\natexlab}[1]{#1}
\providecommand{\url}[1]{\texttt{#1}}
\providecommand{\urlprefix}{URL }
\expandafter\ifx\csname urlstyle\endcsname\relax
  \providecommand{\doi}[1]{\discretionary{}{}{}https://doi.org/#1}\else
  \providecommand{\doi}[1]{\discretionary{}{}{}\urlstyle{rm}\url{https://doi.org/#1}}\fi

\bibitem[{Ponda et~al.(2009)Ponda, Kolacinski, and Frazzoli}]{ponda2009FIMObjectiveFunctions}
Ponda, S., Kolacinski, R., and Frazzoli, E., \enquote{Trajectory optimization for target localization using small unmanned aerial vehicles,} \emph{AIAA guidance, navigation, and control conference}, 2009, p. 6015.

\bibitem[{Xu et~al.(2024)Xu, Zhu, Li, Wu, and Xu}]{xu2024systematical}
Xu, S., Zhu, B., Li, X., Wu, X., and Xu, T., \enquote{Systematical Sensor Path Optimization Solutions for AOA Target Localization Accuracy Improvement with Theoretical Analysis,} \emph{IEEE Transactions on Vehicular Technology}, 2024.

\bibitem[{Xu(2020)}]{xu2020rssAndAOA}
Xu, S., \enquote{Optimal sensor placement for target localization using hybrid RSS, AOA and TOA measurements,} \emph{IEEE Communications Letters}, Vol.~24, No.~9, 2020, pp. 1966--1970.

\bibitem[{Sahu et~al.(2022)Sahu, Wu, Babu, MR, and Ottersten}]{sahu2022optimal}
Sahu, N., Wu, L., Babu, P., MR, B.~S., and Ottersten, B., \enquote{Optimal sensor placement for source localization: A unified ADMM approach,} \emph{IEEE Transactions on Vehicular Technology}, Vol.~71, No.~4, 2022, pp. 4359--4372.

\bibitem[{Shahidian and Soltanizadeh(2017)}]{shahidian2017single}
Shahidian, S. A.~A., and Soltanizadeh, H., \enquote{Single-and multi-UAV trajectory control in RF source localization,} \emph{Arabian Journal for Science and Engineering}, Vol.~42, No.~2, 2017, pp. 459--466.

\bibitem[{Shahidian and Soltanizadeh(2015)}]{shahidian2015autonomous}
Shahidian, S. A.~A., and Soltanizadeh, H., \enquote{Autonomous trajectory control for limited number of aerial platforms in RF source localization,} \emph{2015 3rd RSI International Conference on Robotics and Mechatronics (ICROM)}, IEEE, 2015, pp. 755--760.

\bibitem[{Dogancay(2012)}]{dogancay2012uav}
Dogancay, K., \enquote{UAV path planning for passive emitter localization,} \emph{IEEE Transactions on Aerospace and Electronic systems}, Vol.~48, No.~2, 2012, pp. 1150--1166.

\bibitem[{Hameed et~al.(2022)Hameed, Maqsood, Hashmi, Saeed, and Riaz}]{hameed2022reinforcement}
Hameed, R., Maqsood, A., Hashmi, A., Saeed, M., and Riaz, R., \enquote{Reinforcement learning-based radar-evasive path planning: A comparative analysis,} \emph{The Aeronautical Journal}, Vol. 126, No. 1297, 2022, pp. 547--564.

\bibitem[{Li et~al.(2022)Li, Cheng, and Hu}]{li2022genetic}
Li, W., Cheng, L., and Hu, J., \enquote{Research on stealthy UAV path planning based on improved genetic algorithm,} \emph{2022 International Conference on Artificial Intelligence and Computer Information Technology (AICIT)}, IEEE, 2022, pp. 1--6.

\bibitem[{Zhang et~al.(2020)Zhang, Wu, Dai, and He}]{zhang2020improcedAstar}
Zhang, Z., Wu, J., Dai, J., and He, C., \enquote{Rapid penetration path planning method for stealth uav in complex environment with bb threats,} \emph{International Journal of Aerospace Engineering}, Vol. 2020, No.~1, 2020, p. 8896357.

\bibitem[{Basbous(2018)}]{basbous20182dSimulatedAnnealing}
Basbous, B., \enquote{2D UAV path planning with radar threatening areas using simulated annealing algorithm for event detection,} \emph{2018 international conference on artificial intelligence and data processing (IDAP)}, IEEE, 2018, pp. 1--7.

\bibitem[{Zhao et~al.(2016)Zhao, Niu, Ma, and Ji}]{zhao2016fastAstar}
Zhao, Z., Niu, Y., Ma, Z., and Ji, X., \enquote{A fast stealth trajectory planning algorithm for stealth UAV to fly in multi-radar network,} \emph{2016 IEEE International Conference on Real-time Computing and Robotics (RCAR)}, IEEE, 2016, pp. 549--554.

\bibitem[{Gao et~al.(2014)Gao, Gong, Zhen, Zhao, and Sun}]{gao2014AntColonyOpt}
Gao, C., Gong, H., Zhen, Z., Zhao, Q., and Sun, Y., \enquote{Three dimensions formation flight path planning under radar threatening environment,} \emph{Proceedings of the 33rd Chinese Control Conference}, IEEE, 2014, pp. 1121--1125.

\bibitem[{Kabamba et~al.(2006)Kabamba, Meerkov, and Zeitz~III}]{kabamba2006optimal}
Kabamba, P.~T., Meerkov, S.~M., and Zeitz~III, F.~H., \enquote{Optimal path planning for unmanned combat aerial vehicles to defeat radar tracking,} \emph{Journal of Guidance, Control, and Dynamics}, Vol.~29, No.~2, 2006, pp. 279--288.

\bibitem[{Pelosi et~al.(2012)Pelosi, Kopp, and Brown}]{pelosi2012range}
Pelosi, M., Kopp, C., and Brown, M., \enquote{Range-limited UAV trajectory using terrain masking under radar detection risk,} \emph{Applied Artificial Intelligence}, Vol.~26, No.~8, 2012, pp. 743--759.

\bibitem[{Zabarankin et~al.(2006)Zabarankin, Uryasev, and Murphey}]{zabarankin2006calculusOfVariation}
Zabarankin, M., Uryasev, S., and Murphey, R., \enquote{Aircraft routing under the risk of detection,} \emph{Naval Research Logistics (NRL)}, Vol.~53, No.~8, 2006, pp. 728--747.

\bibitem[{Costley et~al.(2022{\natexlab{a}})Costley, Swedeen, Droge, Christensen, and Leishman}]{costley2022path}
Costley, A., Swedeen, J., Droge, G., Christensen, R., and Leishman, R.~C., \enquote{Path planning with uncertainty for aircraft under threat of detection from ground-based radar,} \emph{The Journal of Defense Modeling and Simulation}, 2022{\natexlab{a}}, p. 15485129251327178.

\bibitem[{Inanc et~al.(2008)Inanc, Muezzinoglu, Misovec, and Murray}]{inanc2008framework}
Inanc, T., Muezzinoglu, M.~K., Misovec, K., and Murray, R.~M., \enquote{Framework for low-observable trajectory generation in presence of multiple radars,} \emph{Journal of guidance, control, and dynamics}, Vol.~31, No.~6, 2008, pp. 1740--1749.

\bibitem[{Sun et~al.(2017)Sun, Liu, Dai, and Grymin}]{sun2017two}
Sun, C., Liu, Y.-C., Dai, R., and Grymin, D., \enquote{Two approaches for path planning of unmanned aerial vehicles with avoidance zones,} \emph{Journal of Guidance, Control, and Dynamics}, Vol.~40, No.~8, 2017, pp. 2076--2083.

\bibitem[{Al-Dahhan and Schmidt(2020)}]{al2020voronoi}
Al-Dahhan, M. R.~H., and Schmidt, K.~W., \enquote{Voronoi boundary visibility for efficient path planning,} \emph{IEEE Access}, Vol.~8, 2020, pp. 134764--134781.

\bibitem[{Judd and McLain(2001)}]{judd2001spline}
Judd, K., and McLain, T., \enquote{Spline based path planning for unmanned air vehicles,} \emph{AIAA Guidance, Navigation, and Control Conference and Exhibit}, 2001, p. 4238.

\bibitem[{Liu et~al.(2022)Liu, Gao, Liu, Liu, and Han}]{liu2022fusion}
Liu, Z., Gao, L., Liu, F., Liu, D., and Han, W., \enquote{Fusion of weighted Voronoi diagram and $\mathrm{A}^{\ast}$ algorithm for mobile robot path planning,} \emph{2022 2nd International Conference on Electrical Engineering and Mechatronics Technology (ICEEMT)}, 2022, pp. 403--406.
\newblock \doi{10.1109/ICEEMT56362.2022.9862795}.

\bibitem[{Choi et~al.(2010)Choi, Curry, and Elkaim}]{choi2010real}
Choi, J.-w., Curry, R., and Elkaim, G., \enquote{Real-time obstacle-avoiding path planning for mobile robots,} \emph{AIAA Guidance, Navigation, and Control Conference}, 2010, p. 8411.

\bibitem[{Zhang et~al.(2018)Zhang, Liu, and Tang}]{zhang2018quantitative}
Zhang, C., Liu, H., and Tang, Y., \enquote{Quantitative evaluation of Voronoi graph search algorithm in UAV path planning,} \emph{2018 IEEE 9th International Conference on Software Engineering and Service Science (ICSESS)}, IEEE, 2018, pp. 563--567.

\bibitem[{Chen et~al.(2016)Chen, Chen, and Xu}]{chen2016path}
Chen, X., Chen, X., and Xu, G., \enquote{The path planning algorithm studying about UAV attacks multiple moving targets based on Voronoi diagram,} \emph{International Journal of Control and Automation}, Vol.~9, No.~1, 2016, pp. 281--292.

\bibitem[{Chen et~al.(2013)Chen, Xiaoqing, Jiyang, and Linfei}]{chen2013research}
Chen, P., Xiaoqing, L., Jiyang, D., and Linfei, Y., \enquote{Research of path planning method based on the improved Voronoi diagram,} \emph{2013 25th Chinese Control and Decision Conference (CCDC)}, IEEE, 2013, pp. 2940--2944.

\bibitem[{Tzoreff and Weiss(2017)}]{tzoreff2017path}
Tzoreff, E., and Weiss, A.~J., \enquote{Path design for best emitter location using two mobile sensors,} \emph{IEEE Transactions on Signal Processing}, Vol.~65, No.~19, 2017, pp. 5249--5261.

\bibitem[{Dai and Cochran~Jr(2010)}]{dai2010path}
Dai, R., and Cochran~Jr, J., \enquote{Path planning and state estimation for unmanned aerial vehicles in hostile environments,} \emph{Journal of guidance, control, and dynamics}, Vol.~33, No.~2, 2010, pp. 595--601.

\bibitem[{Costley et~al.(2022{\natexlab{b}})Costley, Christensen, Leishman, and Droge}]{costley2022sensitivity}
Costley, A., Christensen, R., Leishman, R.~C., and Droge, G.~N., \enquote{Sensitivity of single-pulse radar detection to aircraft pose uncertainties,} \emph{IEEE Transactions on Aerospace and Electronic Systems}, Vol.~59, No.~3, 2022{\natexlab{b}}, pp. 2286--2295.

\bibitem[{Costley et~al.(2023)Costley, Christensen, Leishman, and Droge}]{costley2023sensitivity}
Costley, A., Christensen, R., Leishman, R.~C., and Droge, G., \enquote{Sensitivity of the Probability of Radar Detection to Radar State Uncertainty,} \emph{IEEE Transactions on Aerospace and Electronic Systems}, 2023.

\bibitem[{Li et~al.(2016)Li, Jing, and Bai}]{radarRecognition}
Li, H., Jing, W., and Bai, Y., \enquote{Radar emitter recognition based on deep learning architecture,} \emph{2016 CIE International Conference on Radar (RADAR)}, 2016, pp. 1--5.
\newblock \doi{10.1109/RADAR.2016.8059512}.

\bibitem[{Buccieri et~al.(2009)Buccieri, Perritaz, Mullhaupt, Jiang, and Bonvin}]{Buccieri2009_diffflat}
Buccieri, D., Perritaz, D., Mullhaupt, P., Jiang, Z.-P., and Bonvin, D., \enquote{Velocity-Scheduling Control for a Unicycle Mobile Robot: Theory and Experiments,} \emph{IEEE Transactions on Robotics}, Vol.~25, No.~2, 2009, pp. 451--458.
\newblock \doi{10.1109/TRO.2009.2014494}.

\bibitem[{Boots et~al.(1999)Boots, Okabe, and Sugihara}]{boots1999spatial}
Boots, B., Okabe, A., and Sugihara, K., \enquote{Spatial tessellations,} \emph{Geographical information systems}, Vol.~1, 1999, pp. 503--526.

\bibitem[{Held and de~Lorenzo(2020)}]{held2020efficient}
Held, M., and de~Lorenzo, S., \enquote{An efficient, practical algorithm and implementation for computing multiplicatively weighted Voronoi diagrams,} \emph{arXiv preprint arXiv:2006.14298}, 2020.

\bibitem[{Cox(1972)}]{cox1972numerical}
Cox, M.~G., \enquote{The numerical evaluation of B-splines,} \emph{IMA Journal of Applied mathematics}, Vol.~10, No.~2, 1972, pp. 134--149.

\bibitem[{Bradbury et~al.(2018)Bradbury, Frostig, Hawkins, Johnson, Leary, Maclaurin, Necula, Paszke, Vander{P}las, Wanderman-{M}ilne, and Zhang}]{jax2018github}
Bradbury, J., Frostig, R., Hawkins, P., Johnson, M.~J., Leary, C., Maclaurin, D., Necula, G., Paszke, A., Vander{P}las, J., Wanderman-{M}ilne, S., and Zhang, Q., \enquote{{JAX}: composable transformations of {P}ython+{N}um{P}y programs,} , 2018.
\newblock \urlprefix\url{http://github.com/google/jax}.

\bibitem[{W{\"a}chter and Biegler(2006)}]{IPOPT}
W{\"a}chter, A., and Biegler, L.~T., \enquote{On the implementation of an interior-point filter line-search algorithm for large-scale nonlinear programming,} \emph{Mathematical programming}, Vol. 106, No.~1, 2006, pp. 25--57.

\bibitem[{Norton et~al.(2023)Norton, Stagg, Ward, and Peterson}]{norton2023decentralized}
Norton, T., Stagg, G., Ward, D., and Peterson, C.~K., \enquote{Decentralized sparse gaussian process regression with event-triggered adaptive inducing points,} \emph{Journal of Intelligent \& Robotic Systems}, Vol. 108, No.~4, 2023, p.~72.

\bibitem[{Hart et~al.(1968)Hart, Nilsson, and Raphael}]{hart1968formal}
Hart, P.~E., Nilsson, N.~J., and Raphael, B., \enquote{A formal basis for the heuristic determination of minimum cost paths,} \emph{IEEE transactions on Systems Science and Cybernetics}, Vol.~4, No.~2, 1968, pp. 100--107.

\bibitem[{Dierckx(1995)}]{dierckx1995curve}
Dierckx, P., \emph{Curve and surface fitting with splines}, Oxford University Press, 1995.

\bibitem[{Virtanen et~al.(2020)Virtanen, Gommers, Oliphant, Haberland, Reddy, Cournapeau, Burovski, Peterson, Weckesser, Bright, {van der Walt}, Brett, Wilson, Millman, Mayorov, Nelson, Jones, Kern, Larson, Carey, Polat, Feng, Moore, {VanderPlas}, Laxalde, Perktold, Cimrman, Henriksen, Quintero, Harris, Archibald, Ribeiro, Pedregosa, {van Mulbregt}, and {SciPy 1.0 Contributors}}]{2020SciPy-NMeth}
Virtanen, P., Gommers, R., Oliphant, T.~E., Haberland, M., Reddy, T., Cournapeau, D., Burovski, E., Peterson, P., Weckesser, W., Bright, J., {van der Walt}, S.~J., Brett, M., Wilson, J., Millman, K.~J., Mayorov, N., Nelson, A. R.~J., Jones, E., Kern, R., Larson, E., Carey, C.~J., Polat, {\.I}., Feng, Y., Moore, E.~W., {VanderPlas}, J., Laxalde, D., Perktold, J., Cimrman, R., Henriksen, I., Quintero, E.~A., Harris, C.~R., Archibald, A.~M., Ribeiro, A.~H., Pedregosa, F., {van Mulbregt}, P., and {SciPy 1.0 Contributors}, \enquote{{{SciPy} 1.0: Fundamental Algorithms for Scientific Computing in Python},} \emph{Nature Methods}, Vol.~17, 2020, pp. 261--272.
\newblock \doi{10.1038/s41592-019-0686-2}.

\end{thebibliography}

\end{document}